\documentclass[12pt]{amsart}
\usepackage{geometry}
\usepackage{graphicx}
\usepackage{amssymb}
\usepackage{epstopdf}
\usepackage{amsmath,amscd}
\usepackage{amsthm}
\usepackage{url,verbatim}
\usepackage{stmaryrd}
\usepackage{longtable}
\usepackage[table]{xcolor}

\RequirePackage[colorlinks,citecolor=blue,urlcolor=blue]{hyperref}
\usepackage{breakurl}
\theoremstyle{plain}
\newtheorem{theorem}{Theorem}
\newtheorem{definition}[theorem]{Definition}
\newtheorem{lemma}[theorem]{Lemma}
\newtheorem{proposition}[theorem]{Proposition}
\newtheorem{corollary}[theorem]{Corollary}
\newtheorem{example}[theorem]{Example}

\newtheorem{remark}[theorem]{Remark}
\newtheorem{question}[theorem]{Question}

\newcommand\es{\varnothing}

\newcommand\ol{\overline}
\newcommand\Aut{\mathrm{Aut}}
\newcommand \sA{\mathcal{A}}
\newcommand\sB{{\mathcal B}}
\newcommand\sH{{\mathcal H}}

\newcommand\RR{{\mathbb R}}
\newcommand\ZZ{{\mathbb Z}}

\newcommand\sBp{\sB_{\text{\rm p}}}
\newcommand\Sip{\Si_{\text{\rm b}}}

\newcommand\LL{{\mathbb L}}

\renewcommand\a{\alpha}
\newcommand\om{\omega}

\newcommand\be{\beta}
\newcommand\si{\sigma}
\newcommand\al{\alpha}
\newcommand\eps{\epsilon}

\newcommand\qq{\qquad}
\newcommand\q{\quad}
\newcommand\resp{respectively}
\newcommand\spz{\mathrm{span}}

\newcommand\oo{\infty}
\newcommand\sG{{\mathcal G}}

\newcommand\sQ{{\mathcal Q}}
\newcommand\sK{{\mathcal K}}

\newcommand\sS{{\mathcal S}}

\newcommand\sC{{\mathcal C}}
\newcommand\Ga{\Gamma}
\newcommand\Si{\Sigma}

\newcommand\id{{\bf 0}}

\newcommand\mut{\mu_{\text{\rm T}}}

\newcommand\Sitb{\Si_{\text{\rm T,b}}}

\newcommand\betat{\beta_{\text{\rm T}}}

\renewcommand\ell{l}

\newcommand\pd{\partial}

\newcommand\ghf{graph height function}
\newcommand\hf{height function}
\newcommand\GHF{group height function}

\newcommand\hdi{$\sH$-difference-invariant}

\newcommand\refl{\text{\rm ref}}
\newcommand\rev{\text{\rm rev}}
\newcommand\rr{\text{\rm rr}}
\newcommand\ppm{($\pm$)}

\newcounter{mycount1}\newcounter{mycount2}\newcounter{mycount3}\newcounter{mycount}
\newenvironment{romlist}{\begin{list}{\rm(\roman{mycount1})}%
   {\usecounter{mycount1}\labelwidth=1cm\itemsep 0pt}}{\end{list}}

\newenvironment{numlist}{\begin{list}{\rm\arabic{mycount2}.}%
   {\usecounter{mycount2}\labelwidth=1cm\itemsep 0pt}}{\end{list}}
\newenvironment{letlist}{\begin{list}{\rm(\alph{mycount3})}%
   {\usecounter{mycount3}\labelwidth=1cm\itemsep 0pt}}{\end{list}}

\numberwithin{equation}{section}
\numberwithin{theorem}{section}
\numberwithin{figure}{section}
\numberwithin{table}{section}

\title{SAWs, trails, kites, patterns}
\author{Geoffrey R.\ Grimmett}
\address{Statistical Laboratory, Centre for
Mathematical Sciences, Cambridge University, Wilberforce Road,
Cambridge CB3 0WB, UK} 
\email{grg@statslab.cam.ac.uk}
\urladdr{\url{http://www.statslab.cam.ac.uk/~grg/}}

\begin{document}

\begin{abstract}
A \emph{self-avoiding walk} (SAW) on a vertex-transitive graph $G$ 
is a path that visits no vertex more than once; a \emph{trail}
is a path that traverses no edge more than once.  The number of $n$-step self-avoiding walks (\resp, trails)
starting at a given vertex has order $\mu(G)^n$ (\resp, $\mut(G)^n$), where $\mu(G)$
is called the connective constant and $\mut(G)$ is the trail constant. It is immediate that $\mu=\mut$ if $G$ either has degree $3$ or less, or is a tree. We prove a necessary and sufficient condition for the
strict inequality $\mu < \mut$ if $G$ has degree $4$ or more and is not a tree, 
and in addition $G$ possesses a so-called height function. This condition is satisfied whenever the
height function is transitive.
The core of the proof is a pattern theorem for graphs possessing height functions, and this is achieved
using an unpublished argument for $\ZZ^d$ by Hammersley, reported by Janse van Rensburg et al.\
(J.\ Stat.\ Phys.\ 85, 1996, 103--130). Some ancillary results are included concerning SAWs on the 
line graph of $G$.
\end{abstract}

\date{October 14, 2026} 

\keywords{Self-avoiding walk, trail, connective constant, vertex-transitive graph, height function, 
pattern theorem, line graph}
\subjclass[2010]{05C30, 60K35, 82B20}
\maketitle

\section{Introduction}\label{sec:int}

A \emph{self-avoiding walk} (abbreviated to SAW) on a graph $G=(V,E)$ is a path that visits no vertex
more than once. 
SAWs were introduced in the  theory of polymerization (see Orr \cite{Orr} and Flory \cite{f}),
and their critical behaviour has attracted the 
abundant attention since of mathematicians and physicists 
(see, for example, the book of Madras and Slade \cite{ms}
and the lecture notes \cite{bdgs}). 

Let $\sG$ denote the set of all infinite, connected,
locally finite, (vertex-)transitive graphs without loops or parallel edges.
For $G\in\sG$ we designate some vertex the \emph{root} of $G$, denoted $\id$. 
Let $\si_n$
be the number of $n$-step SAWs on $G$ starting  at $\id$.
It is a fundamental theorem of Hammersley \cite{jmhII} that the limit
\begin{equation*}
\mu=\mu(G) := \lim_{n \to \oo} \si_n^{1/n}
\end{equation*}
exists, and the limit $\mu$ is called the \emph{connective constant} of $G$. 
(The assumption  of transitivity can be weakened to quasi-transitivity, \cite{jmhII}.)
It is proved in \cite{GL-Comb} that
\begin{equation*}
\sqrt{d-1} \le \mu(G) \le d-1,
\end{equation*}
where $d$ is the vertex-degree of $G$.

SAWs are vertex-avoiding in that they visit no vertex twice or more. In contrast, a \emph{trail} 
(sometimes called a \lq self-avoiding trail') on $G\in \sG$ is a 
path that is edge-avoiding in that it traverses no edge (irrespective of direction) twice or more.  Let $\tau_n$ be the number of $n$-step trails
starting at the root. The study of the sequence  $(\tau_n)$ was initiated by Malakis \cite{Mal1,Mal2},
who noted the existence  of the \lq trail constant'
\begin{equation*}
\mut=\mut(G):= \lim_{n\to\oo} \tau_n^{1/n}.
\end{equation*}
Since every SAW is a trail, we have $\mu \le \mut$.
This provokes the question of determining for which graphs the strict inequality $\mu<\mut$ holds.
Guttmann \cite{AJG85,AGII} pointed out, \emph{inter alia}, that $\mu=\mut$ if $G$ has degree $3$, since in that case (with the 
possible exceptions of the endpoints of the walk)  a walk is vertex-avoiding if and only if it is edge-avoiding.
It is immediate that $\mu=\mut$ if $G$ is a  tree, and thus we are led to the following question.

\begin{question}
Let $G\in \sG$. Is it the case that $\mu <\mut$ if and only if $G$ has degree $4$ or more and possesses a cycle?
\end{question}

The main theorem of the current work is a partial resolution of this question.

\begin{theorem}\label{thm:trail}
Suppose $G\in\sG$ has  a transitive \ghf\ $(h,\sH)$. If $G$ has degree $4$ or more and possesses 
a cycle, then $\mu<\mut$.
\end{theorem}

The meaning of \lq graph height function' is outlined in the next paragraph, but we mention 
in advance that the forthcoming Theorem \ref{prop:3} provides a (verifiable) necessary
and sufficient condition for a graph $G\in\sG$ with a \ghf\ (either transitive or not)
to satisfy $\mu<\mut$.

A \lq \ghf' (abbreviated henceforth as simply \lq height function') of $G$
is a function $h$ that maps the vertex-set of $G$ onto $\ZZ$ with  certain properties that
will be explained in Section \ref{sec:patt2}. In brief,
a \hf\ is a formalization of the notion of the \lq heights' of vertices, and it permits a SAW analysis 
based around so-called bridges.
Height functions
have been very useful in demonstrating properties of SAWs and connective constants on general vertex-transitive graphs.
They were introduced in \cite{GL-loc}, and for further details the reader is referred to the 
review \cite{GL-amen}. (For an elaboration of the assumption in Theorem \ref{thm:trail}
that the \hf\ be transitive
rather than merely quasi-transitive, 
see Section \ref{sec:trailthm} of the current work.)

Possession of a \hf\ seems the norm for the graphs in $\sG$
that are mostly considered  in the literature.  For example, it is shown in \cite{GL-loc} that
$\ZZ^d$, the regular trees, the hexagonal lattice (and indeed the triangular lattice)
have \hf s.
It seems one has to work quite hard to obtain a graph in $\sG$ with no \hf, and we note
in Section \ref{ssec:crystals} 
that the Cayley graphs of the Grigorchuk and Higman graphs are such (see \cite{GL2017}).

In principle, the proof of Theorem \ref{thm:trail}
allows a strictly positive lower bound for $\mut-\mu$, but this is unlikely to be
useful in practice. By numerical methods, 
one obtains that $\mut(\ZZ^2) \approx 2.72$
and $\mut(\ZZ^3) \approx 4.84$ (see, for example, \cite{CG93,SZW}), 
whereas $\mu(\ZZ^2)\approx 2.64$ and $\mu(\ZZ^3) \approx 4.68$
(see, for example, \cite{Clisby,AGII,JSG} and also \cite[Sect.\ 7.3]{Hug}).
Furthermore, Guttmann \cite{AGII} outlined a proof of the strict inequality
$\mu(\ZZ^d) < \mut(\ZZ^d)$ for $d \ge 2$; this proof is completed 
and extended in the current paper.  

We outline next our approach to the proof of Theorem \ref{thm:trail}. Suppose $G\in\sG$  has degree $4$ and
possesses a cycle. A trail on $G$ may enter a cycle, circumnavigate this cycle, and exit at the same point as it arrived.
Such a local configuration constitutes a special type of \lq pattern'. 
The pattern theorem of Kesten \cite{HK63}, albeit
restricted to SAWs on the hypercubic lattices $\ZZ^d$, suggests that such patterns should exist in all but
a small proportion of $n$-step trails, and furthermore the latter proportion should be exponentially small in $n$.
This is enough to imply $\mu(G)<\mut(G)$.

In this paper we present a pattern theorem for trails in a graph that possesses a \hf. The proof does not follow
the classic proof of Kesten \cite{HK63}, but is instead motivated by an unpublished argument of Hammersley \cite{Ham96},
as adapted in another setting by Janse van Rensburg et al.\ \cite{Rens96}. Hammersley's approach is in harmony
with the existence  of a \hf\ since it is based on a renewal theory of bridges noted by Kesten \cite{HK63}.

In Section \ref{sec:ham}, we present an account of what we believe may be not too distant from
 Hammersley's idea for patterns in
$\ZZ^d$; this is included here both for historical reasons and also since it may be useful in other settings. 
Section \ref{sec:patt2} contains an extension of
this new proof to a substantially more general family of graphs than the hypercubic lattices, 
namely vertex-transitive graphs possessing a \hf.
The SAW pattern theorem of that section has a  cousin theorem for trails, proved similarly. 
The proof of the general pattern theorem is given in Section \ref{sec:proofpatt1}.
It then remains to show in Section \ref{sec:trailthm}
that, for the graphs under study, 
all but an exponentially small proportion of trails contain a pattern that cannot belong to a 
SAW; the pattern in question is the circumnavigation of an entire cycle.

The pattern theorem of Section \ref{sec:patt2} (Theorem \ref{thm:pattern1})
differs from that of Section \ref{sec:ham} (Theorem \ref{thm:pattern0}) in
that it applies to a much larger class of graphs than \lq merely' the hypercubic lattices $\ZZ^d$.
A further distinction is that, whereas Theorem \ref{thm:pattern0} is proved using  both \emph{reflection} and \emph{reversal} of walks, Theorem \ref{thm:pattern1} is proved using \emph{reversal} only. 

There is a substantial literature devoted to strict inequalities for critical points. Much of that work
lies in the context of models of probability theory, for which \lq finite-energy-type' methods may often 
be used
(see, for example, \cite{AG}).
The situation with SAWs and trails is rather different because of the 
unforgiving combinatorial nature  of their definitions.

There is some relationship between the results of this paper and those of \cite{GL14}. The latter contains strict
inequalities between SAW
connective constants subject to perturbations in the underlying graph, by either adding new 
edges, or factoring  by a normal subgroup of the automorphism group. Paper \cite{GL14} proceeds by an algebraic 
adaptation of the pattern theorem, in the classical style of Kesten's proof from \cite{HK63,ms}. 
In contrast, a principal
purpose of the current paper is to highlight Hammersley's alternative approach 
(using bridges) to the pattern theorem.

Let $G\in\sG$ have vertex-degree $d$. The line graph (or
\lq covering' graph) of $G$ may be viewed as the graph obtained by
replacing each vertex $v$ of $G$ by a copy of the complete graph $K_d$,
each vertex of which is the mid-point of one of the edges incident to $v$. 
The line graph of $G$ is denoted $G_K$ in Section \ref{sec:line}.
When $d=3$, $G_K$ is the so-called \lq Fisher graph' of $G$, as
obtained by the \lq Fisher transformation' of \cite{fisher}.
This transformation is a valuable tool in the study of a number of
interactive models on graphs (see, for example, \cite{GL13, MW}).

In an alternative approach to studying trails, one observes that a trail on $G\in\sG$
gives rise to a certain SAW on its line graph $G_K$. 
When $d\ge 4$, not all SAWs on $G_K$ correspond to trails on $G$.
In Section \ref{sec:line}, we include certain ancillary results concerning SAWs on line graphs,
including a relation between the connective constants of $G$ and $G_K$.

The methods of this paper may be applied to obtain strict inequalities for
connective constants arising in certain other SAW-type problems. Here are two examples.
\begin{letlist}
\item 
Let $k\ge 1$.
Walks on $G$ that visit no vertex more than $k$ times (sometimes called
 $k$-\emph{tolerant} walks, see \cite[Sect.\ 7.4.4]{Hug}).
\item
Let $k,l\ge 1$.
Walks $w=(w_0,w_1,\dots)$ on $G$ such that $d_G(w_i,w_j) \ge k$ whenever
$|i-j|\ge l$. Here, $d_G$ denotes graph-distance.
\end{letlist}

\section{Preliminaries}\label{sec:not}

A graph $G=(V,E)$ has automorphism group $\Aut(G)$, and it is 
called \emph{vertex-transitive} (abbreviated to \emph{transitive}) if, for all $u,v\in V$, 
there exists $\al\in\Aut(G)$  such that $\al(u)=v$. It is called \emph{quasi-transitive} if
$\Aut(G)$ acts on $V$ with only finitely many orbits, and \emph{locally finite}
if each vertex is incident to only finitely many edges.
We denote by $\sG$ the set of all infinite, connected, locally finite, transitive graphs without loops or parallel edges.

For $G=(V,E)\in\sG$, we call some vertex of $V$ the \emph{root} and denote it $\id$.
An edge of $E$ with endvertices $u$, $v$ is written $\langle u,v\rangle$. 
If $\langle u,v\rangle \in E$,
vertices $u$ and $v$ are called \emph{neighbours}, and we write $u \sim v$.
The \emph{neighbourhood} $\pd u$ of a vertex $u$ is
the set of its neighbours.
Since our graphs are undirected,
we identify $\langle u,v\rangle$ and $\langle v,u \rangle$. 
The \emph{degree} $\deg(u)$ of vertex $u$ is the number $|\pd u|$ of its neighbours.
Since $G$ is assumed transitive, $\deg(u)$ is constant for $u \in V$, and its common value is
called the \emph{degree} of $G$. The term $\sG_d$ denotes the set of all $G\in\sG$ with degree $d$. 

A \emph{walk} on $G\in\sG$ is an alternating sequence 
$w=(w_0,e_0,w_1,e_1,\dots)$ of vertices $w_i$ and edges $e_i
=\langle w_i, w_{i+1}\rangle$; since $G$ has no parallel edges, we lose no information by
 suppressing mention of the edges $e_i$ in the above listing of $w$. 
 A walk $w$ with $n$ edges is said to have $n$ \emph{steps}, and we denote this
 number by $\vert w\vert$.
 The walk $w$ is called \emph{self-avoiding} if $w_i\ne w_j$ for all $i\ne j$,
and is called a \emph{trail} if $e_i \ne e_j$ for all $i\ne j$.  
The walk $w=(w_0,w_1,\dots,w_n)$ is called \emph{closed} if $w_0=w_n$, and
this closed walk is a \emph{cycle} if $(w_0,w_1,\dots,w_{n-1})$ is
self-avoiding. The notation $\ZZ_k$ denotes a cycle of length $k$, that is,
the Cayley graph of the cyclic group $\ZZ_k$.

Let $\si_n$ be the number of $n$-step self-avoiding walks (or \lq  SAWs')
from $\id$. 
It is standard that $\si_{m+n} \le \si_m\si_n$, whence the limit
$\mu(G) = \lim \si_n^{1/n}$ exists, and is called the 
\emph{connective constant}.

Let $d \ge 2$. For $u\in\ZZ^d$ we write $u=(u_1,u_2,\dots,u_d)$.
The \emph{hypercubic lattice} $\ZZ^d$ has vertex-set $V=\ZZ^d$ and edge-set 
$E$ containing all edges $\langle u,v\rangle$
with $\sum_{i=1}^d |u_i-v_i| = 1$. The \emph{root} (or origin) of $\ZZ^d$ is the zero-vector 
$\id=(0,0,\dots,0)$.
 The reader is referred to \cite{ms} for an account  of SAWs on $\ZZ^d$,
and to \cite{GL-amen} for more general graphs.

When working on $\ZZ^d$, we define $h:\ZZ^d\to \ZZ$ by $h(x)=x_1$, the first coordinate of $x\in\ZZ^d$.
Assume $n\ge 1$ in the following.
A \emph{half-space walk} $w=(w_0,w_1,\dots,w_n)$ on $\ZZ^d$ is a walk such that 
$h(w_i) > h(w_0)$ for all $i\ge 1$.
The SAW $w=(w_0,w_1,\dots,w_n)$ on $\ZZ^d$ is called a \emph{bridge} 
if $w_0=\id$ and $0< h(w_i)\le h(w_n)$ for $1\le i\le n$.
It is standard that the number $b_n$ of $n$-step bridges on $\ZZ^d$ satisfies $b_{m+n} \ge b_mb_n$,
and hence  the limit
$\be(\ZZ^d) := \lim b_n^{1/n}$ exists, and is called the \emph{bridge constant}.
Furthermore,
\begin{equation}\label{eq:0}
b_n \le \be(\ZZ^d)^n, \qq n\ge 1.
\end{equation}

We shall encounter some more general definitions in Section \ref{sec:patt2}.
For real numbers $x$, $y$ we write
\begin{equation*}
x\wedge y = \min\{x,y\}, \qq x\vee y = \max\{x,y\}.
\end{equation*}
When using the symbols $\pm$, $\mp$, we adopt the usual convention that one sticks
throughout the relevant sentence with either the upper symbol or the lower symbol.

\section{A pattern theorem for hypercubic lattices}\label{sec:ham}

Let $d\ge 2$, and consider the graph $\ZZ^d$ embedded in the usual way in $\RR^d$.
For $x\in\RR$ we define the hyperplane $H_x=\{x\}\times\RR^{d-1}$.
Bridges of $\ZZ^d$ were defined in Section \ref{sec:not}.

\begin{definition}\label{def:cut}
Let $e_k=\langle w_k,w_{k+1}\rangle$ be an edge of the walk $w=(w_0,w_1,\dots,w_n)$
satisfying $h(w_k)\ne h(w_{k+1})$,
and write $m_k= \frac12(h(w_k)+h(w_{k+1}))$.
The edge $e_k$ is called a \emph{cut-edge} of $w$ if  $e_k$ is the unique edge of $w$ that intersects $H_{m_k}$.
A bridge is called \emph{prime} if no edge except the first is a cut-edge.
\end{definition} 

Let $\sB$ (\resp, $\sBp$) be the set of all bridges (\resp, prime bridges).
A bridge may be decomposed into prime bridges, as explained next. 
Let $w=(w_0=\id,w_1,\dots,w_n)$ 
be a bridge and let $e_k=\langle w_k,w_{k+1}\rangle$ be the 
earliest cut-edge of $w$ other than $e_0$  (if no such cut-edge exists, we take $k=n$).
Let $\al_1\in\sBp$ be the sub-SAW from $w_0$ to $w_k$, noting that $\al_1$
is a prime bridge. Let $w'$
be the bridge obtained from the sub-SAW of $w$ from $w_k$ to $w_n$ by translating it
in such a way that its first vertex $w_k$ 
is translated to the root; that is $w'=(w_k,w_{k+1}\dots,w_n)-w_k$.  
We may now express $w$ in the form  $w=\al_1\ast w'$.
By iteration, we find that $w=\al_1\ast\al_2\ast\cdots\ast \al_r$ for 
prime bridges $\al_i\in\sBp$;
this is called the \emph{prime decomposition} of $w$ (see \cite[Sect.\ 4]{HK63} and \cite[Sect.\ 4.2]{ms}). 
Here, $w$ is reconstructed from the $\al_i$ by a sequential augmentation:
one starts with $\al_1$, and iteratively for each $k$, one augments the path $\pi_{k-1}$ obtained 
from $\al_1,\al_2,\dots,\al_{k-1}$ by adding to it the SAW $\al_k+v_{k-1}$,
where $v_{k-1}$ is the endpoint of $\pi_{k-1}$ other than the root $\id$. 

We turn to so-called patterns.
A SAW $p=(p_0,p_1,\dots,p_r)$, where $r\ge 1$,  is called a \emph{pattern}
if if has no cut-edge (this definition may be relaxed slightly but we retain it for later convenience).
The pattern $p$ is said to occur at the point $w_i$ of the SAW $w=(w_0,w_1,\dots,w_n)$ if
$i+r\le n$ and 
$(w_i,w_{i+1},\dots,w_{i+r}) = w_i+p-p_0$.
The \emph{number of occurrences of} $p$
in $w$ is the cardinality of the set of $i$ such that $p$ occurs at $w_i$.
A pattern is called \emph{bridgeable} if there exists some bridge $w$ such that $p$ occurs at least once in $w$.
Let $\Sip$ denote the set of all bridgeable patterns.

The arguments of this section make use of the operations of reflection and reversal of SAWs.
Let $\rho:\ZZ^d\to \ZZ^d$ be the \emph{reflection mapping} given by
\begin{align*}
\rho(x_1,x_2,\dots,x_d)=(-x_1,x_2,\dots,x_d).
\end{align*} 
For a SAW $w=(w_0, w_1,\dots,w_n)$ we define the reflected SAW $w_\refl=\rho(w)$ and the reversed SAW $w_\rev = R(w)$ by
\begin{align}
w_\refl &:= \rho(w) = (\rho(w_0),\rho(w_1),\dots,\rho(w_n)),\label{eq:refl}\\
w_\rev &:=R(w)=(w_n,w_{n-1},\dots,w_0).\label{eq:rev}
\end{align}
Note that $\rho$ and $R$ commute, and write
\begin{equation}\label{eq:rr}
w_\rr:= R(\rho(w)) = \rho(R(w)). 
\end{equation}

It may be checked by the definition of cut-edge that, for a SAW $p$,
\begin{equation}\label{eq:allornone}
\text{either all or none of $p$, $p_\refl$, $p_\rev$, $p_\rr$ are patterns.}
\end{equation}

We note an important fact.

\begin{lemma}\label{lem:oneb}
Let $w$ be a bridge. Any pattern $p$ may occur in
at most one of the prime bridges 
in the prime decomposition of $w$.
\end{lemma}

\begin{proof}
This holds since such $p$ has no cut-edge.
\end{proof}

\begin{theorem}[Pattern theorem]\label{thm:pattern0}
Let $d \ge 2$ and let $p\in\Sip$. Let $\si_n(P,\eps)$ be the number
of $n$-step SAWs starting at $\id$ that 
contain at most $\eps n$ copies of any of the elements of the set $P=\{p, p_\refl, p_\rev, p_\rr\}$. 
There exists
$a>0$ such that 
\begin{equation}\label{eq:100}
\limsup_{n\to\oo}  \si_n(P,\eps)^{1/n} < \mu(G),
\end{equation}
whenever $0\le \eps<a$.
\end{theorem}

The reason for including  the patterns $p_\refl$, $p_\rev$, $p_\rr$ in the set $P$ of \eqref{eq:100}
is as follows. As in \cite{HW62}, the proof has three stages (presented here in reverse order). 
\begin{letlist}
\item
Using a geometrical construction, a SAW $w$ is transformed
into two half-space walks $h_1$, $h_2$. In doing so, part of $w$ is reversed.

\item
Each half-space walk $h_i$ is transformed into a bridge denoted $c_i$. 
This is done by decomposing $h_i$ into bridges and \lq negative bridges', 
and then reflecting the negative bridges to obtain \lq ordinary' bridges.

\item
The central part of the proof is to show that the number of $n$-step bridges containing no copy of $p$
is exponentially small compared with the total number of $n$-step bridges.

\end{letlist}
In order that there be no copies of $p$ in the bridges $c_i$, one needs to rule out copies
of both $p$ and $p_\refl$ in the half-space walks $h_i$. Working back through stage (a), 
one requires that $w$ contain no copies of $p$, $p_\refl$, $p_\rev$, $p_\rr$.

\begin{remark}\label{rem:7}
Whereas Theorem \ref{thm:pattern0} makes use of both reflection and reversal, the more general
Theorem \ref{thm:pattern1} of Section \ref{sec:patt2} makes use of reversal but not reflection.
When studying $\ZZ^d$ in this section, we choose  
to utilise the bridge arguments of Hammersley and Welsh \cite{HW62}
 rather than  the more general approach of Grimmett and Li
\cite{GL-loc}. In summary, the conclusion of the
above theorem holds with the references to $p_\refl$ and $p_\rr$
removed throughout.
\end{remark}

\begin{remark}[Relationship to Kesten's pattern theorem]\label{rem:8}
Kesten's theorem, as presented in \cite{ms}, applies to copies of a SAW $p$ of $\ZZ^d$ for which there exists a
SAW containing at least three copies of $p$ (see \cite[Prop.\ 7.1.3]{ms}). 
For graphs with a \hf, a pattern satisfies the last condition if and only if it is bridgeable.  
Kesten's theorem
makes no assumption on the cut-edges of $P$, whereas Theorem \ref{thm:pattern0} requires there to be none,
and additionally the latter bounds the number of appearances of members of the set  
$\{p, p_\refl, p_\rev, p_\rr\}$. As remarked in Remark \ref{rem:7},
the reflected patterns
$p_\refl$, $p_\rr$ may be omitted here. The principal virtues of the current theorem 
compared to that of \cite{HK63, ms} lie in the relative simplicity and clarity
of its proof, and its extension (given in the next section) to a class of graphs
much larger than just $\ZZ^d$.
\end{remark}

\begin{remark}[History]\label{rem:1}
The first pattern theorem was proved by Kesten in his influential work \cite{HK63},
in which he introduced also the use of renewal theory. The proof given there
is somewhat complex, and may be found essentially in \cite[Chap.\ 3]{ms}. That theorem has proved 
fundamental in the study of SAWs; see, for example, \cite{Abu,GL14,Rens96,ms}. 
Madras \cite{Mad99} has studied related but distinct combinatorial problems which require
a different proof of the corresponding pattern theorem. 

The work \cite{Rens96} of Janse van Rensburg et al.\ is directed at self-avoiding ribbons rather than SAWs, and it
uses a proof of the pattern theorem of which the main idea is attributed by the authors to Hammersley
\cite{Ham96}. 
It was in reading \cite{Rens96} that the current author learned of Hammersley's argument.
The core of the argument appears also in work \cite{AGII} of Guttmann from 1985, with
a credit given to one of the authors of \cite{Rens96}. 
\end{remark}

Let $p\in\Sip$, and let $b_n(p,0) $ (\resp, $b_n(P,0)$) 
be the number of $n$-step bridges that contain no copies of $p$ (\resp, 
copies of none of the elements of the set $P=\{p, p_\refl, p_\rev, p_\rr\}$); 
we adopt a similar notation ($\si_n(p,0)$, etc)
in replacing bridges by SAWs. 
We have in the usual way (using Lemma \ref{lem:oneb}) that
\begin{align}\label{eq:116}
\si_{m+n}(P,0) &\le \si_{m}(P,0)\si_{n}(P,0),\\
b_{m+n}(P,0) &\ge b_{m}(P,0)b_{n}(P,0), 
\end{align}
so that, by the subadditive limit theorem, the following limits exists:
\begin{align}\label{eq:106}
 \mu(P) &:= \lim_{n\to\oo} \si_n(P,0)^{1/n},\\
 \beta(P) &:= \lim_{n\to\oo} b_n(P,0)^{1/n},\label{eq:106a}\\
  \beta(p) &:= \lim_{n\to\oo} b_n(p,0)^{1/n}.\label{eq:107}
\end{align} 
By the definitions, we have
\begin{equation}\label{eq:115}
\mu(P) \le \mu(p), \qq \beta(P) \le \beta(p).
\end{equation}
 
We state and prove an important proposition which is at the heart of Hammersley's approach. 

\begin{proposition}\label{lem:1}
Let $d \ge 2$ and let $p\in\Sip$. 
We have the strict inequality
\begin{equation}\label{eq:110}
\beta(p) < \be(\ZZ^d),
\end{equation}
and, moreover,
\begin{equation}\label{eq:108}
\beta(P)= \mu(P),\qq \beta(\ZZ^d)=\mu(\ZZ^d).
\end{equation}
\end{proposition}

\begin{proof} 
The basis of the argument is a renewal-theoretic formulation which originated in Kesten \cite{HK63}.
As noted above, a bridge $w$ may be expressed uniquely as a concatenation of its successive prime elements, namely 
$w=\al_1\ast\al_2\ast\cdots\ast \al_r$.
Let $b_n$ be the number of $n$-step bridges and $q_n$ the number of $n$-step,  prime bridges.
As in classical renewal theory
(see, for example, \cite[Chap.\ 10]{GS20}, \cite[Sect.\ 4]{HK63}, \cite[Sect.\ 4.2]{ms}, \cite{Rens96}), we have that
\begin{equation}\label{eq:102}
b_n=q_n+\sum_{m=1}^{n-1} q_m b_{n-m}, \qq n\ge 1.
\end{equation}
Therefore, the generating functions
\begin{equation*}
B(x)=\sum_{\si\in\sB} x^{|\si|},\qq Q(x)=\sum_{\si\in\sBp} x^{|\si|},
\end{equation*}
satisfy
\begin{equation*}
B(x) = Q(x) + Q(x)B(x),
\end{equation*}
whence
\begin{equation}\label{eq:ren}
B(x) = \frac{Q(x)}{1-Q(x)}\qq \text{whenever } Q(x)<1.
\end{equation}

By the definition of the bridge constant $\be=\be(\ZZ^d)$, $B$ has radius of convergence $1/\be$, 
so that,
by Abel's lemma (see, for example, \cite[p.\ 171]{GS20}),  
\begin{equation}\label{eq:101}
Q(\be^{-1})=1.
\end{equation}
This equation was first derived by Kesten \cite[Thm 5]{HK63}.

Let $p\in\Sip$, and let $Q_p$ be the generating function of prime bridges containing no copy
of $p$. 
Since $p$ is assumed bridgeable, by \eqref{eq:101},
$Q_p (\be^{-1}) <1$.
Since $Q_p$ is continuous inside its ball of convergence, there exists $\eta=\eta(p)>0$ such that 
\begin{equation}\label{eq:103}
Q_p(\be^{-1}+\eta)<1.
\end{equation}

We repeat the above for bridges containing no copies of
$p$.  By the argument leading to \eqref{eq:ren} 
adapted to such bridges, we deduce by \eqref{eq:103} that $\be(p)^{-1}>\be^{-1}+\eta$, as required for \eqref{eq:110}.
That is, the proportion of $n$-step bridges that contain no copies of $p$ is exponentially small in $n$.
That completes the first part of this proof. The remaining
part is fairly standard and the details are omitted in favour of a summary.

The second equality of \eqref{eq:108} is 
well known (and was first proved by Hammersley and Welsh \cite{HW62}). 
The reader is referred to the overview of their proof presented after the statement of Theorem \ref{thm:pattern0}.
The proof details are omitted here, but may be found at \cite[Cor.\ 3.1.8]{ms},
\cite[Sec.\ 2]{bdgs}, \cite[p.\ 436]{Hug}, and elsewhere; moreover, this argument is cited in \cite[Lemma 3.1]{Rens96}.
Using Lemma \ref{lem:oneb} and the invariance of $P$ under reflection and reversal,
it may be checked that the same argument may be applied to SAWs that contain no copy of any 
element of the set $P=\{p, p_\refl, p_\rev, p_\rr\}$.

As in \cite[Thm 3.1.1, eqn 3.1.8]{ms}, we conclude that there exists $c=c(p,d)<\oo$ such that
\begin{equation}\label{eq:105}
\si_n(P, 0) \le e^{c\sqrt n} \beta(P)^n
\end{equation}
for all large $n$. 

It is immediate that $\beta(P)\le \mu(P)$, 
and the first equality of \eqref{eq:108} follows by \eqref{eq:105}.
\end{proof}

\begin{proof}[Proof of Theorem \ref{thm:pattern0}]
By Proposition \ref{lem:1} and \eqref{eq:115}, we have 
\begin{equation*}
\mu(P) = \be(P) \le \be(p) < \mu(\ZZ^d),
\end{equation*}
whence
\begin{equation*}
\lim_{n\to\oo} \si_n(P,0)^{1/n} < \mu(\ZZ^d).
\end{equation*}
The extended claim \eqref{eq:100} follows
in the same manner as \cite[Lemma 7.2.5]{ms},
to which readers desirous of a full proof are referred.
\end{proof}

\section{A pattern theorem for general crystals}\label{sec:patt2}

The word \lq crystal' in the section title is taken from \cite{GG26a,Ham61}.

\subsection{Height functions and bridges}\label{ssec:crystals}

Height functions were introduced under the title \lq\ghf' in \cite{GL-loc}.
Since the distinction between \lq\ghf' and \lq\GHF' is irrelevant here,
we adopt the briefer term \lq height function'. 
First we define a \hf, and then we use such functions to define bridges in general graphs.

\begin{definition} \cite[Defn 3.1]{GL-loc}\label{def:height}
Let $G =(V,E)\in \sG$ have  root labelled $\id$. 
\begin{romlist}
\item
A \emph{\hf} on $G$ is a pair $(h,\sH)$ such that:
\begin{letlist}
\item $h:V \to\ZZ$, and $h(\id)=0$, 
\item $\sH\le \Aut(G)$ acts quasi-transitively on $G$ 
such that $h$ is \emph{\hdi}, in the sense that
$$
h(\a v) - h(\a u) = h(v) - h(u), \qq \a \in \sH,\ u,v \in V,
$$
\item for  $v\in V$,
there exist $u,w \in \pd v$ such that
$h(u) < h(v) < h(w)$.
\end{letlist}
\item
A \hf\ $(h,\sH)$ is called \emph{transitive} (\resp, \emph{unimodular}) if the action
of $\sH$ is transitive (\resp, unimodular).
\end{romlist}
\end{definition}

Note that $(h,\sH)$ is a \hf\ if and only if $(-h,\sH)$ is a \hf. We shall sometimes omit explicit reference to the group $\sH$.

Height functions have proved useful for several problems associated with SAWs, since they allow
arguments based around bridges. For example, it was shown in \cite[Thm 5.1]{GL-loc} that a graph possessing
a unimodular \hf\
has connective constant with the property of locality (the condition of unimodularity was removed 
by the main theorem of \cite{Lin20}). 

The existence, non-existence, and consequences  of \hf s 
were explored for various families of
graphs in \cite{GG26a,GL2017,GL-loc,GL20}; see \cite{GL-amen} for a review.
See \cite[Chap.\ 8]{LyP} and \cite[eqn (3.1)]{GL-loc} 
for discussions of unimodularity. For clarity, we list some examples of graphs with or without a \hf.
\begin{letlist}
\item
The hypercubic lattice $\ZZ^d$ and the triangular lattice have transitive \hf s (\cite{GL-loc}).
\item
The hexagonal lattice has a quasi-transitive \hf\ but no transitive \hf\ (\cite{GL-loc}).
\item
The Cayley graphs of the first Grigorchuk group of \cite{Grig80}
and the Higman group of \cite{Hg} have no \hf\  
(\cite[Thms 5.1, 8.1]{GL2017}).
\end{letlist}

Henceforth in this section, let $G\in\sG$ and let $(h,\sH)$ be a \hf\ for $G$.
A walk $w=(w_0,w_1,\dots)$ is called \emph{increasing} (\resp, \emph{non-decreasing})
if $h(w_i)<h(w_{i+1})$ for $i\ge 0$ (\resp, $h(w_i)\le h(w_{i+1})$ for $i\ge 0$);
similar definitions hold for decreasing SAWs.
The \emph{span} of a SAW $w=(w_0,w_1,\dots,w_n)$ on $G$ is given by
\begin{equation}\label{eq:span}
\spz(w) =\max_i h(w_i) - \min_i h(w_i). 
\end{equation}

We turn to bridges.
There are two types of bridge depending on whether we work 
with the \hf\ $h$ or the \hf\ $-h$.
A SAW $w=(w_0=\id,w_1,\dots,w_n)$ with $n\ge 1$
is called a ($+$)\emph{bridge} if 
\begin{equation}\label{eq:b+}
h(w_0)< h(w_i) \le h(w_n), \qq 0< i \le n,
\end{equation}
and a ($-$)\emph{bridge} if
\begin{equation}\label{eq:b-}
h(w_0)> h(w_i) \ge h(w_n), \qq 0< i \le n.
\end{equation}
 
\begin{definition}\label{def:cut2}
Let $e_k=\langle w_k,w_{k+1}\rangle$ be an edge of the walk $w=(w_0,w_1,\dots,w_n)$.
The edge $e_k$ is called a  
\begin{align*}
\text{{\rm($+$)}\emph{cut-edge} of $w$ if}\q  &h(w_i)
\begin{cases} \le h(w_k) &\text{for $i\le k$},\\
> h(w_k) &\text{for $i>k$},
\end{cases}\\
\text{{\rm($-$)}\emph{cut-edge} of $w$ if}\q  &h(w_i)
\begin{cases} \ge h(w_k) &\text{for $i\le k$},\\
< h(w_k) &\text{for $i>k$},
\end{cases}
\end{align*}
A ($\pm$)bridge $w=(w_0,w_1,\dots,w_n)$ is called \emph{prime} if no edge 
except its first is a ($\pm$)cut-edge.
\end{definition} 

Let $b_n^\pm$ be the number of $n$-step ($\pm$)bridges.
We write $\sB^\pm$ (\resp, $\sBp^\pm$) for the set of ($\pm$)bridges 
(\resp, prime ($\pm$)bridges). 

It may be shown \lq in the usual way' (see \cite[Sect.\ 4]{GL-loc} and Lemma \ref{lem:10}) 
that the limits $\beta^\pm=\lim_{n\to\oo} (b_n^\pm)^{1/n}$
exist, and $\beta^\pm =\be^\pm(G)$ is called the ($\pm$)\emph{bridge constant} of $G$. 
Moreover,  it is proved in \cite[Prop.\ 4.2]{GL-loc} that there exists $r=r(G,h,\sH)<\oo$ such that
\begin{equation}\label{eq:bbnd}
b_n^\pm  \le \be^\pm(G)^{n+r}, \qq n \ge 1.
\end{equation}
If the \hf\ $(h,\sH)$ is transitive, we may take $r=0$; we shall explain the value of $r$ 
for the non-transitive case in Section \ref{ssec:qtrans}.
It is proved in \cite[Thm 4.3]{GL-loc} that $\beta^+=\beta^-$ if the \hf\ is unimodular. 
In addition the common value is 
independent of the choice of \hf\ $(h,\sH)$, and moreover it equals $\mu(G)$; 
this is generally false when $G$ is non-unimodular (see \cite{Lin20}). 

The following was proved by Lindorfer \cite{Lin20} using 
an extension of the methods of \cite{GL-loc,HW62}. The set $\sQ$ denotes the
set of members of $\sG$ enriched by relaxing the assumption of transitivity to quasi-transitivity.

\begin{theorem}[\mbox{\cite[Thm 3.2]{Lin20}}]\label{thm:bridge}
Let $G\in \sQ$ possess a  \hf\ $(h,\sH)$. Then $\mu(G)=\beta^+(G) \vee\beta^-(G)$.
\end{theorem}

Outlines of the proof of this theorem when the \hf\ is transitive (\resp, quasi-transitive)
are provided within the proofs of Section \ref{sec:proofpatt1}.

\begin{figure}
\centerline{\includegraphics[width=0.5\textwidth]{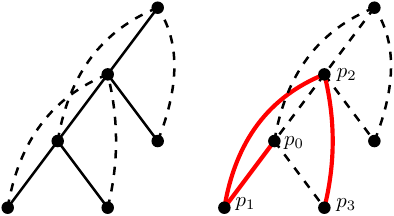}}
\caption{The grandparent graph of Trofimov \cite{trof} is derived from the binary tree by adding an edge from each vertex 
to its grandparent. Part of this graph is illustrated, with the pattern $p$ indicated on the right.}
\label{fig:grand}
\end{figure}

A \emph{pattern} is defined as a SAW of length one or more that has neither a ($+$)cut-edge
nor a ($-$)cut-edge.
A pattern $p$ is ($\pm$)\emph{bridgeable} if
there exists a ($\pm$)bridge $w$ and $\al\in\sH$ such that $w$ contains $\al (p)$.
Let $\Sip^\pm$ denote the set of ($\pm$)bridgeable patterns, and $\Sip=\Sip^+\cup \Sip^-$.
One may wonder whether $\sH$ could be replaced by $\Aut(G)$ in the above definition
of bridgeability,
but we note that bridges without images of $p$ under $\Aut(G)$ are rarer than those without 
images under $\sH$.

\begin{example}\label{ex:nonuni}
Let $p$ be a pattern.
When $G=\ZZ^d$, it is not hard to see that $p\in\Sip^+$ if and only if $p\in \Sip^-$.
This does not hold for all $G$, as the following example indicates. Consider the grandparent graph $G$
of Figure \ref{fig:grand}.  Endow $G$ with the horocyclic height function
of \cite[Rem.\ 3.3]{GL-loc}, and let $p$ be the highlighted SAW. Then $p\in\Sip^-$
but $p\notin \Sip^+$.
\end{example}

\subsection{Pattern theorem for SAWs}\label{ssec:psaw}

Unlike the hypercubic lattices, general graphs do not usually have a property of reflection-invariance.
The action of reflection is then replaced (as in \cite{GL-loc}) by that of reversal.
As in \eqref{eq:rev}, the \emph{reversal} of a SAW $w=(w_0,w_1,\dots,w_n)$ is the SAW
\begin{equation*}
w_\rev:=(w_n,w_{n-1},\dots,w_0),
\end{equation*}
and we write $w_\rev = R(w)$. 

We shall employ the notation of \eqref{eq:106}--\eqref{eq:107}, and we recall \eqref{eq:115}.
There follows the pattern theorem for graphs with a \hf. It concerns general bridgeable patterns 
except those lying in a 
certain special set $\Pi$ of SAWs satisfying
\begin{equation}\label{eq:es}
\text{$\Pi=\es$ when the \hf\ is transitive.}
\end{equation}
When the \hf\ is not transitive but quasi-transitive, the excluded set $\Pi$ is a certain set of SAWs given
in Definition \ref{def:pi}. For definiteness, we note that
\begin{equation}\label{eq:pibnd}
0 \le |\Pi| \le M(M-1),
\end{equation}
where $M$ is the number of orbits of $V$ under $\sH$. The set $\Pi$ of excluded patterns
is not given uniquely in terms of $G$ and the \hf, and thus, in considering a specific pattern $p$, 
it may be possible to engineer $\Pi$ to exclude $p$. See Section \ref{ssec:qtrans} and particularly 
Remark \ref{rem:pi}.

We write
\begin{equation}\label{eq:401}
\mu(p\cup p_\rev) = \limsup_{n\to\oo} \sigma_n(p\cup p_\rev,0)^{1/n},
\end{equation}
noting that the $\limsup$ is a limit as in \eqref{eq:106} when the \hf\ is transitive.
The last turns out to hold in the quasi-transitive case also so long as $p\notin \Pi$, see Corollary \ref{cor:1}.

\begin{theorem}[Pattern theorem for SAWs]\label{thm:pattern1}
Let $(h,\sH)$ be a \hf\ of $G\in\sG$,
and let $p\in\Sip\setminus\Pi$. 
\begin{letlist}
\item 
We have that $\mu(p\cup p_\rev,0) < \mu(G)$.

\item

Let $\si_n(p\cup p_\rev,\eps)$ be the number of $n$-step SAWs starting at $\id$ that contain 
no more than $\eps n$ copies of either $p$ or $p_\rev$. There exists
$a>0$ such that 
\begin{equation}\label{eq:400}
\limsup_{n\to\oo} \si_n(p\cup p_\rev,\eps)^{1/n}< \mu(G),
\end{equation}
whenever $0\le\eps<a$.
\end{letlist}
\end{theorem}

See Section \ref{sec:proofpatt1} for the proof.

\subsection{Pattern theorem for trails}

Key to the proof of Theorem \ref{thm:trail} 
is the observation that the above pattern theorem for SAWs applies, 
\emph{mutatis mutandis}, to trails. When replacing SAWs by trails, the notions of height function, bridge,
half-space walk/trail, reversal/reflection, cut-edge, primality, bridgeability, and so on, 
carry over with no significant change. In order
to distinguish between notation for SAWs and for trails, in the latter case we adorn relevant notation with the suffix \lq T\rq.

Let $G\in\sG$ have \hf\ $(h,\sH)$.
A \ppm\emph{t-bridge} (or \lq trail-bridge') is a trail $w=(w_0,w_1,\dots,w_n)$ satisfying 
\eqref{eq:b+} or \eqref{eq:b-}, as appropriate.
A \ppm t-bridge is \emph{prime} if
it has no ($\pm$)cut-edge apart from its first edge.
A \emph{t-pattern} (or \lq trail-pattern')  is a trail with length $1$ or more containing neither
($+$)cut-edge nor ($-$)cut-edge.
A t-pattern $p$ is \ppm\emph{t-bridgeable} if there exists a \ppm t-bridge $w$ containing some image of $p$
under $\sH$.
The set of ($\pm$)t-bridgeable t-patterns is denoted $\Sitb^\pm$, and we let $\Sitb=\Sitb^+\cup\Sitb^-$.
 
We note the trail analogue of Theorem \ref{thm:bridge}; the proof is similar and is omitted.

\begin{theorem}\label{thm:t-bridge}
Let $G\in \sG$ possess a  \hf\ $(h,\sH)$. Then $\mut(G)=\betat^+(G) \vee\betat^-(G)$.
\end{theorem}

As a consequence of this theorem, $\betat^+\vee\betat^-$ is invariant
with respect to the choice of height function.
Next is the pattern theorem for trails.
The set $\Pi$, viewed as a set of bridgeable trails, is as in Theorem \ref{thm:pattern1}, and satisfies 
\eqref{eq:es} as there.

For $p\in\Sitb$, we write
\begin{equation}\label{eq:401t}
\mut(p\cup p_\rev) = \limsup_{n\to\oo} \tau_n(p\cup p_\rev,0)^{1/n},
\end{equation}
noting that the $\limsup$ is a limit as in \eqref{eq:106} when the \hf\ is transitive.
As in the case of SAWs (Corollary \ref{cor:1}),
the $\limsup$ is indeed a limit in the quasi-transitive case whenever $p\notin \Pi$
(where $\Pi$ is as in Section \ref{ssec:psaw}).

\begin{theorem}[Pattern theorem for trails]\label{thm:trail2}
Let $(h,\sH)$ be a \hf\ of $G\in\sG$,
and let $p\in\Sitb\setminus\Pi$. 
\begin{letlist}
\item 
We have that 
$\mut(p\cup p_\rev,0) < \mut(G)$.

\item

Let $\tau_n(p\cup p_\rev,\eps)$ be the number of $n$-step trails starting at $\id$ that contain 
no more than $\eps n$ copies of either $p$ or $p_\rev$. There exists
$a>0$ such that 
\begin{equation}\label{eq:450}
\limsup_{n\to\oo} \tau_n(p\cup p_\rev,\eps)^{1/n}< \mut(G),
\end{equation}
whenever $0\le \eps<a$.
\end{letlist}
\end{theorem}

\begin{figure}
\centerline{\includegraphics[width=0.4\textwidth]{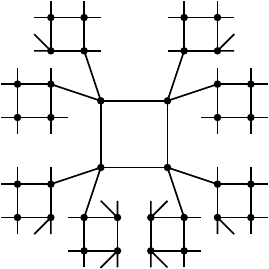}}
\caption{Part of the free product graph $G=\ZZ_4*K_2*K_2$. Any trail that circumnavigates the central square must pass along two siblings.}
\label{fig:zkk}
\end{figure}

The proofs of Theorems \ref{thm:t-bridge} and \ref{thm:trail2}  follow largely those of
Theorems \ref{thm:bridge} and \ref{thm:pattern1} and are omitted. See also Section \ref{ssec:comments}.

\begin{example}\label{rem:no}
One naturally asks whether such a theorem holds for graphs without a \hf. The obvious proof route
would go via Kesten's approach, a simple interpretation of which might be phrased as follows.
``For any given t-pattern $p$ and any sufficiently long trail
$w$, one may perturb $w$ locally to create a copy of $p$.'' This assertion
is problematic even for a graph with a quasi-transitive \hf.

\begin{figure}
\centerline{\includegraphics[width=0.3\textwidth]{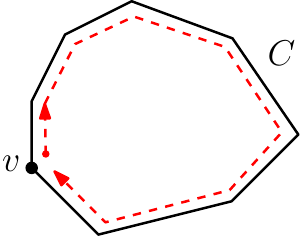}}
\caption{With $C$ a cycle and $v\in C$, the t-pattern $p(C,v)$
is a circumnavigation of $C$ starting and ending at $v$.}
\label{fig:pdelta}
\end{figure}

Consider the 
graph $G=\ZZ_m*K_2*K_2$ of Figure \ref{fig:zkk} with, for concreteness, $m=4$;
 here $\ZZ_m$ denotes an $m$-cycle and $K_2$ an edge. 
(See \cite{Gilch} for basic properties of free product graphs.)
Let $\square$ denote a given $4$-cycle of $G$ and let $v$ be a vertex 
of $\square$. The t-pattern $p(\square,v)$ is defined as the trail that circumnavigates 
$\square$ from $v$ to $v$ (the direction of travel is immaterial). See Figure \ref{fig:pdelta}

An edge of $G$ is called \emph{cyclic} if it lies in some cycle.
Two non-cyclic edges of $G$ are called \emph{siblings} if they have a common vertex.
Let $w$ be a SAW that uses no pair of siblings, and whose first and last edges are cyclic. There exists no trail 
$w'$ with the same endpoints as $w$ that contains a copy (under $\Aut(G)$) of $p(\square,v)$. On the other hand
there exists a trail containing three disjoint copies of $p(\square,v)$; cf.\ Remark \ref{rem:8}.

It may be seen that $G$ possesses a \hf\ with respect to which $p(\square, v)$ is ($+$)bridgeable.
Indeed, as in the proof of the forthcoming Theorem \ref{prop:3},
for every \hf\ of $G$ there exists $v\in \square$ such that $p(\square,v)$ is ($+$)bridgeable.
\end{example}

\section{Proof of Theorem \ref{thm:pattern1}}\label{sec:proofpatt1}

For clarity of exposition, we shall first assume that
the \hf\ is transitive (that is, $\sH$ acts transitively on the vertex-set $V$), 
and then in Sections \ref{ssec:qtrans}--\ref{ssec:qtpf} 
we shall explain the further steps in the non-transitive case. This
mimics the approach of \cite{GL-loc,Lin20}.

Let $G\in\sG$ have \hf\ $(h,\sH)$, and let
\begin{equation}\label{eq:defd1}
d=d(h)=\max\bigl\{|h(u)-h(v)|: u,v\in V,\ u \sim v\bigr\}.
\end{equation}

\subsection{The proof in the case of transitive \hf}
\label{ssec:trans}

We assume in this subsection that the height function is transitive, and
we recall from \eqref{eq:es} that $\Pi=\es$ in this case.

\begin{lemma}\label{lem:24}
The pattern $p$ satisfies $p\in\Sip^\pm$ if
and only if $p_\rev\in \Sip^\mp$. 
\end{lemma}

\begin{proof}
We may assume without loss of generality that  $p \in\Sip^+$, and
we shall deduce that $p_\rev=R(p)$ is ($-$)bridgeable.
Let $w=(w_0=\id,w_1,\dots,w_n)$ be a ($+$)bridge containing a copy of $p$. 
By adding a further edge if needed,
we may assume that $h(w_n)> h(w_i)$ for $0\le i<n$.
By the assumed transitivity of the \hf, there exists 
$\al\in\sH$ such that $\al(w_n)=\id$. Then $R(\al w)$ is a ($-$)bridge containing $p_\rev$, as required.
\end{proof}

Let $b_n^\pm$ be the number of $n$-step ($\pm$)bridges,
and, for $p\in\Sip^\pm$, let $b_n^\pm(p,0)$ (\resp, $b_n^\pm(p\cup p_\rev)$)
be the number of such bridges that contain no copy of $p$ (\resp, no copy of either $p$ or $p_\rev$).
Recall $\mu(p\cup p_\rev)$ from \eqref{eq:401}.
Furthermore, let
\begin{align}\label{eq:321}
\beta^\pm(p\cup p_\rev) = \lim_{n\to\oo} b_n^\pm(p\cup p_\rev,0)^{1/n},\qq
\beta^\pm(p) = \lim_{n\to\oo} b_n^\pm(p,0)^{1/n}.
\end{align}
By the transitivity of the \hf, these limits exist as in \eqref{eq:106a}--\eqref{eq:107}.

\begin{proposition}\label{lem:21}\mbox{\hfil}
\begin{letlist}
\item
If $p\in\Sip^\pm$, then
\begin{equation}\label{eq:210}
\beta^\pm(p\cup p_\rev) \le \beta^\pm(p) < \beta^\pm.
\end{equation}

\item
If $p\in\Sip$, then
\begin{equation}\label{eq:208}
\mu(p\cup p_\rev) < \mu(G).
\end{equation}
\end{letlist}
\end{proposition}

\begin{proof}
We follow the proof of Proposition \ref{lem:1}. The proof is given only in outline, and the reader is
referred to \cite{GL-loc, Lin20} for further details.

(a) Inequality \eqref{eq:210} holds in the same manner as  \eqref{eq:110}, and the details
are omitted.

(b)
For \eqref{eq:208} a new argument is needed, and this is provided by the reversal approach of \cite{GL-loc}
as adapted in \cite{Lin20} and illustrated in Figures \ref{fig:half} and \ref{fig:full}. 
We summarise first how the argument of \cite{Lin20} leads
to the equality 
\begin{equation}\label{eq:903}
\mu=\be^+\vee\be^-
\end{equation}
 of Theorem \ref{thm:bridge}. Evidently,
$\mu\ge \be^+\vee \be^-$, and it suffices for \eqref{eq:903} to prove that
\begin{equation}\label{eq:211}
\mu\le \be^+\vee \be^-.
\end{equation}

\begin{figure}
\centerline{\includegraphics[width=0.5\textwidth]{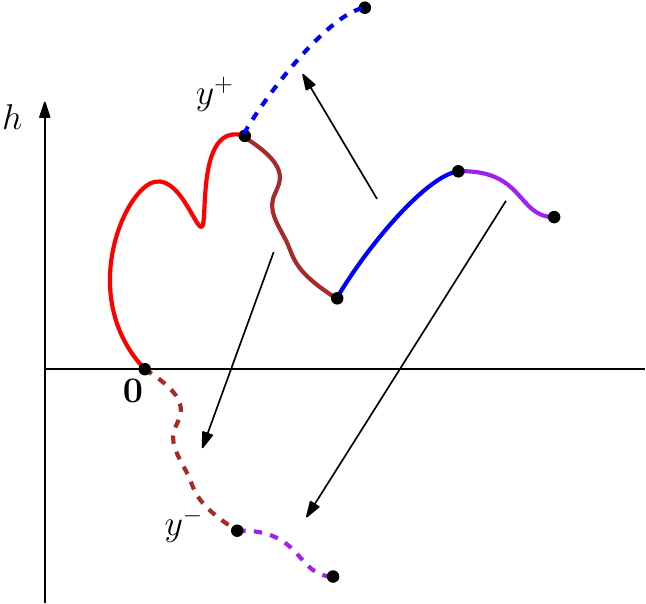}}
\caption{An illustration of the transformation from a half-space walk
to a $(+)$bridge $y^+$ and a $(-)$bridge $y^-$.}\label{fig:half}
\end{figure}

Let $v=(v_0,v_1,\dots,v_m)$ and $w=(w_0,w_1,\dots,w_n)$ be walks on $G$.
Since $\sH$ acts transitively, there exists $\alpha\in\sH$ such that $\al(w_0)=v_m$.
By $v\circ w$ we denote the concatenation of $v$ followed by 
$(\alpha(w_1),\alpha(w_2),\dots,\alpha(w_n))$. 

A \emph{half-space walk} $w=(w_0,w_1,\dots,w_n)$ on $G$ is a walk such 
$h(w_i) > h(w_0)$ for all $i\ge 1$.
Let $H_n$ denote the set of $n$-step half-space SAWs starting at $\id$, and let $w\in H_n$. 
The idea is to decompose $w$ as an alternating sequence of
($+$)bridges $x_1,x_3,\dots$ and ($-$)bridges $x_2,x_4,\dots$, whose spans are decreasing,
and to reposition these to obtain
a ($+$)bridge $y^+:= x_1\circ x_3\circ\cdots$ and a ($-$)bridge $y^-:= x_2\circ x_4\circ\cdots$. 
This is illustrated in Figure \ref{fig:half}.

 Let $a_i$ denote the span of the ($+/-$)bridge $x_i$.
Each such $w$ is specified  by knowledge of $y^+$, $y^-$, and the sequence $(a_i)$.
Let $h_n=|H_n|$ and let $h_n(a_1,a_2,\dots,a_k)$ be the number of elements of $H_n$ with $a_i=\spz(x_i)$
for $1\le i \le k$.
As in \cite[Lemma 3.3]{Lin20},
\begin{align}\label{eq:212}
h_n 
= \sum_{k=1}^n \sum_{\substack{a_1>a_2>\cdots>a_k>0\\ a_1+a_2+\cdots+ a_k\le dn}}
h_n(a_1,a_2,\dots,a_k).
\end{align}
On using the transformation illustrated in Figure \ref{fig:half}, together with the supermultiplicativity of
bridge counts,
\begin{align}\label{eq:2121}
h_n&\le  \sum_{k=1}^n \sum_{\substack{a_1>a_2>\cdots>a_k>0\\ a_1+a_2+\cdots+ a_k\le dn}}
 \sum_{m=0}^{n} b_m^+(a_1+a_3+\cdots) b_{n-m}^-(a_2+a_4+\cdots)\\
&\le P(dn) \sum_{k=1}^n \sum_{m=0}^n  b_m^+b_{n-m}^-, \nonumber
\end{align}
where $b_s^\pm(r)$ denotes the number of $s$-step ($\pm$)bridges with span $r$,
and $P(K)$ is the number of (ordered) partitions of the integer $K$.
Using the Hardy--Ramanujan partition theorem of \cite{HR}, we deduce that there exist constants $A,B>0$ such that,
by \eqref{eq:bbnd}, 
\begin{equation}\label{eq:213}
h_n \le e^{\pi\sqrt{dn/3}}\sum_{m=0}^n b_m^+ b_{n-m}^-
\le A e^{B\sqrt n} (\be^+\vee\be^-)^n.
\end{equation}

\begin{figure}
\centerline{\includegraphics[width=0.5\textwidth]{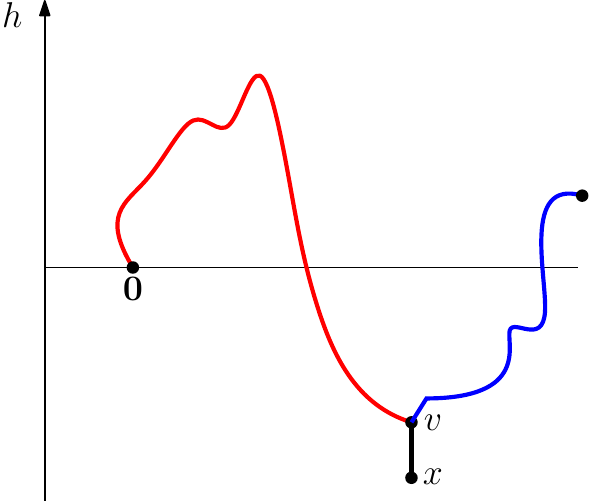}}
\caption{A SAW from $\id$ gives rise to two half-space walks starting at the additional vertex $x$.
The leftmost such walk is a reversal of the original.}\label{fig:full}
\end{figure}

Next we relate half- and full-space SAWs, as in \cite{HW62,Lin20}, for example. 
See Figure \ref{fig:full}.
Let $w$ be a SAW of $G$ starting at $\id$, 
and find its last vertex $v$
with minimal height. By breaking $w$ at $v$ into two sub-walks (the first of which is reversed)
and adding a new edge $\langle v,x\rangle$ to an additional vertex $x$ with strictly smaller height
than $v$, we find that
\begin{equation}\label{eq:214}
\si_n \le \sum_{m=0}^n h_{n-m} h_{m+1},
\end{equation}
whence \eqref{eq:211} follows by \eqref{eq:213}. The equality of \eqref{eq:903} follows.

We next repeat the above argument for walks containing copies of neither $p$ nor $p_\rev$.
We may assume without loss of generality
that $p\in \Sip^+$; in the other case when $p\in\Sip^-$, we recall Lemma \ref{lem:24}
and work instead with $p_\rev$.
In the usual notation, inequalities \eqref{eq:213} and \eqref{eq:214} become (using \eqref{eq:bbnd})
\begin{align}\label{eq:216}
h_n(p\cup p_\rev,0) 
&\le e^{\pi\sqrt{dn/3}}\sum_{m=0}^n b_m^+(p\cup p_\rev,0) b_{n-m}^-(p\cup p_\rev,0)\\
&\le A e^{B\sqrt n} \bigl[\be^+(p\cup p_\rev)\vee\be^-(p\cup p_\rev)\bigr]^n,\nonumber\\
\si_n(p\cup p_\rev,0) 
&\le \sum_{m=0}^n h_{n-m}(p\cup p_\rev,0) h_{m+1}(p\cup p_\rev,0).\label{eq:215}
\end{align}
Therefore, by \eqref{eq:210}, Lemma \ref{lem:24}, and Theorem \ref{thm:bridge},
\begin{align}\label{eq:new3}
\mu(p\cup p_\rev) &\le  \be^+(p\cup p_\rev)\vee \be^-(p\cup p_\rev)\\
&\le  \be^+(p)\vee \be^-(p_\rev)\nonumber\\
&< \be^+\vee \be^- = \mu(G),\nonumber
\end{align}
as claimed in \eqref{eq:208}.
\end{proof}

Part (a) of Theorem \ref{thm:pattern1} is now proved, and part (b) follows as in the proof of Theorem \ref{thm:pattern0},
that is, by the argument of \cite[Lemma 7.2.5]{ms}.

\subsection{The quasi-transitive case}
\label{ssec:qtrans}

Let $G=(V,E)\in\sG$ have \hf\ $(h,\sH)$.
We assume in this and the next subsection that $\sH$ acts quasi-transitively but not transitively, 
and we begin by defining the set $\Pi$ of banned patterns. We shall see that $\Pi$ is not defined uniquely
in terms of $G$ and its height function.

The novel difficulty in the quasi-transitive case, and its solution, 
are exemplified as follows. Let $v=(v_0,v_1,\dots,v_m)$
and $w=(w_0,w_1,\dots,w_n)$ be walks with $h(v_0)<h(v_m)$
and $h(w_0)<h(w_n)$, say. We may wish to extend $v$ by appending an
element of $\sH w$, while maintaining height
differences in a natural way. In the transitive case, one may choose $\al\in\sH$ such that $\al(w_0)=v_m$, and  the 
extended path $v \circ w$ is obtained by adjoining $v$ and $\al(w)$. This path has the property that 
\begin{equation*}
h(v_0)<h(v_m)=h(\al(w_0))<h(\al(w_n)).
\end{equation*}
In the quasi-transitive case, such $\al\in\sH$ may not exist. Instead, we find $\al\in\sH$ such that
$\al(w_0)$ is close to $v_m$, and in addition $h(\al(w_0))>h(v_m)$; then we
find a path from $v_m$ to $\al(w_0)$ using only vertices with heights
strictly between $h(v_m)$ and $h(\al(w_0))$ (apart from its endvertices). Since $V$ has only finitely many
orbits under $\sH$, the lengths of such paths may be taken to be uniformly bounded.
This is formalized as follows.

Let $G/\sH=(\ol V, \ol E)$ be the (simple) quotient graph with vertex-set the set of orbits of $V$ under $\sH$, and
with edge-set as follows. For $v\in V$, write $\ol v=\sH v$. For $v,w\in V$, we place an edge between
$\ol v$ and $\ol w$ if and only if  some member of $\ol w$ is joined by an edge of $G$ to $v$
(if $w\in\ol v$, this edge is a loop).
See \cite[Sect.\ 3.4]{GL14} and \cite[Sect.\ 2]{GL-loc} for discussions of this quotient graph.
Since $G/\sH$ is connected, it has a spanning tree $T$.
Let $L=L(T)$ be a lift of $T$ into $G$, so that $L$ is a sub-tree of $G$ with $M$ vertices,
each of which represents a given orbit.

\begin{figure}
\centerline{\includegraphics[width=0.25\textwidth]{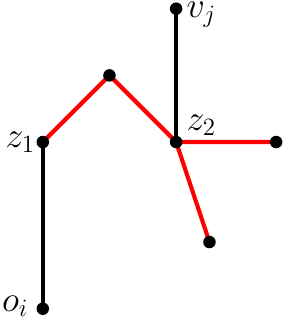}}
\caption{An illustration of the SAW $\nu(o_i,o_j)$ in Proposition \ref{prop:pi}.
The two black lines are SAWs that are increasing (upwards), 
and the red lines form a copy $\alpha(L)$ of the lifted tree $L$.}
\label{fig:nu}
\end{figure}

\begin{proposition}\label{prop:pi}
Let $o_1,o_2,\dots,o_M$ be representatives of the orbits of $V$ under $\sH$.
\begin{letlist}
\item
For $i\ne j$, there
exists a vertex $v_j \in \sH o_j$ 
such that $h(o_i)<h(v_j)$, together with  a SAW, denoted $\nu(o_i,v_j)$, from $o_i$
to $v_j$  all of whose vertices $x$, other than
its endvertices, satisfy $h(o_i)<h(x)< h(v_j)$. We set $\nu(o_i,o_i)$ to
be the $0$-step SAW $\{o_i\}$. 

\item
There exists  $\alpha\in\sH$ such that the following holds.
The SAW $\nu(o_i,o_j)$
comprises three linked edge-disjoint 
sub-SAWs. The first is an increasing SAW from $o_i$ to some vertex $z_1\in \alpha(L)$;
the second is a sub-SAW of $\alpha(L)$ from $z_1$ to some vertex $z_2$; the third is an increasing
SAW from $z_2$ to $v_j$.
\end{letlist}
\end{proposition}

\begin{proof}
See \cite[Prop.\ 3.2]{GL-loc} and its proof.
\end{proof}

The construction  of Proposition \ref{prop:pi}
is illustrated in Figure \ref{fig:nu}. The SAWs $\nu(o_i,o_j)$ and 
their images under $\sH$ are called \emph{connectors}. For $i\ne j$, $\nu(o_i,o_j)$ is a ($+$)bridge.

Associated with a \hf\ $(h,\sH)$
are two integers $d$, $r$. The constant $d$ is as in \eqref{eq:defd1}.
If $\sH$ acts transitively, we set
$r=0$. Otherwise, 
\begin{equation}\label{eq:defr}
r :=\max\bigl\{\vert \nu(o_i,o_j)\vert : i,j\in\{1,2,\dots,M\}\bigr\}.
\end{equation}

The idea of the main proof is to reproduce the geometrical constructions of Section \ref{ssec:trans} 
with the difference that, when re-positioning a sub-SAW, we may need to include a suitable connector
at the join.  
Such connectors are ($+$)bridges which may themselves 
contain copies of a given pattern $p$ (or its reversal $p_\rev$),
and this may disturb the overall pattern count. 
We choose therefore to exclude from consideration patterns that lie in any connector.

\begin{definition}\label{def:pi}
Denote by $\Pi$ the set of patterns $p\in\Sip$ such that either $p$ or $p_\rev$ occurs in some connector.
\end{definition}

With $\Pi$ given thus,
\begin{gather}\label{eq:sym}
\text{every $p\in\Pi$ is a SAW,}\\
\text{$p\in \Pi$ if and only if $p_\rev\in\Pi$.}
\label{eq:sym2}
\end{gather} 
The set $\Pi$ depends on $G$, $(h,\sH)$, and moreover on the choices of $T$, $L$, and the paths $\nu(o_i,o_j)$. 

\begin{remark}\label{rem:pi}
With a given set $\sC$ of connectors, there is an associated set $\Pi=\Pi(\sC)$ of banned patterns, that is,
patterns that are excluded from consideration in the pattern theorem, Theorem \ref{thm:pattern1}. 
Proposition \ref{prop:pi} indicates one way to construct a set $\sC$ of connectors, but there are others.
When a particular pattern $p$ is of interest, one may seek connectors $\sC$  for which $p,p_\rev\notin\Pi(\sC)$.
The situation is especially simple if $\Pi(\sC)=\es$, and this occurs if one may choose every 
$\nu(o_i,o_j)$ to be a non-decreasing SAW. 
This property, when true, may often be checked by inspection of $G$ and the given \hf.
\end{remark}

\begin{figure}
\centerline{\includegraphics[width=0.6\textwidth]{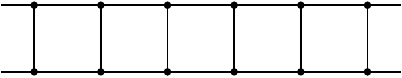}}
\caption{The ladder graph.}
\label{fig:ladd}
\end{figure}

\begin{example}[The ladder graph]\label{ex:ladder}
The ladder $G=\LL$ of Figure \ref{fig:ladd} provides a simple illustration of the above discussion.
Let $\LL=(V,E)$ have vertex-set $V=\{0,1\}\times\ZZ$ and edge-set the set of all pairs of vertices separated by unit distance.
Clearly $\LL$ is transitive, with automorphism group generated by the translation $\tau$ and reflection $\rho$
given by
\begin{equation*}
\tau(v)= (v_1,v_2+1),\q \rho(v) = (1-v_1,v_2), \qq \text{for } v=(v_1,v_2)\in V.
\end{equation*}
Let $h(v) = v_2$ and let $\sH$ be the infinite cyclic group generated by $\tau$, so that $(h,\sH)$
is a quasi-transitive (but not transitive) \hf. 
The lift $L$ of the spanning tree of $G/\sH$ may be taken to be the \lq vertical' edge $e$ between $(0,0)$ and
$(0,1)$. By Remark \ref{rem:pi}, we have that $\Pi=\es$.
\end{example}

\subsection{The proof in the quasi-transitive case}\label{ssec:qtpf}

We next outline the extra steps
needed beyond the arguments of Section \ref{ssec:trans}. 
Further details may be found in the relevant parts of
\cite{GL-loc,Lin20}.

We introduce a slightly more general version of a bridge, to be called an \lq ebridge'. 
A ($\pm$)\emph{ebridge} is a SAW $w=(w_0,w_1,\dots,w_n)$ satisfying
either \eqref{eq:b+} or \eqref{eq:b-} as appropriate; we do not require ebridges to start at $\id$. 
An $(i,j)$-ebridge is a ($+$)ebridge that starts at $o_i$ and ends in $\sH o_j$.
Let $b_n^{i,j}$ be the number of $n$-step $(i,j)$-ebridges,
and $b_n^{i,\cdot} = \sum_j b_n^{i,j}$. 
An $(i,i)$-ebridge $w=(w_0,w_1,\dots,w_n)$ 
is called $i$-\emph{prime} if it is minimal in the sense that it contains no sub-SAW 
$(w_k,w_{k+1},\dots,w_n)$ with $k\ge 1$
that is an $(i,i)$-ebridge. 

Recall the number $M$ of orbits of $V$ under $\sH$.
As before,
$b_n^\pm$ denotes the number of $n$-step \ppm bridges from $\id$,
and $b_n^\pm(p,0)$ (\resp, $b_n^\pm(p\cup p_\rev,0)$)
the number of such bridges that contain no copy of the pattern $p$ (\resp, no copy of either $p$ or $p_\rev$).
Recall \eqref{eq:401} and \eqref{eq:321}, and  write
\begin{gather}\label{eq:323}
\mu(p\cup p_\rev) = \limsup_{n\to\oo} \sigma_n(p\cup p_\rev,0)^{1/n},\\
\beta^\pm(p\cup p_\rev) = \limsup_{n\to\oo} b_n^\pm(p\cup p_\rev,0)^{1/n},\q
\beta^\pm(p) = \limsup_{n\to\oo} b_n^\pm(p,0)^{1/n}.
\label{eq:324}
\end{gather}
We shall see in Lemma \ref{lem:11}(c) that, for $p\notin\Pi$,  the limits superior in
\eqref{eq:324} are actually limits. By Corollary \ref{cor:1},
the 
$\limsup$ in \eqref{eq:323} is a limit when $p\notin\Pi$.

\begin{lemma}\label{lem:10} 
The limit $\be^\pm=\lim_{n\to\oo} (b_n^\pm)^{1/n}$ exists.
Furthermore, we have for $i,j\in\{1,2,\dots,M\}$ that 
\begin{equation}\label{eq:603}
\lim_{n\to\oo} (b_n^{i,\cdot})^{1/n} \to  \beta^+,
\qq \limsup_{n\to\oo}(b_n^{i,j})^{1/n} = \be^+.
\end{equation}
\end{lemma}

\begin{lemma}\label{lem:11}
\mbox{\hfil}
\begin{letlist}
\item
The pattern $p$ satisfies $p\in\Sip^+$ if and only if $p_\rev\in\Sip^-$.

\item
Let $p\in\Sip^+$.
For $i\in\{1,2,\dots,M\}$ and some $n$,
there exists an $n$-step $(i,i)$-ebridge containing a copy of $p$.

\item
If $p\in\Sip^+\setminus\Pi$ and $i\in\{1,2,\dots,M\}$, the limits
\begin{equation*}
\beta^{i,\cdot}(p\cup p_\rev) := \lim_{n\to\oo} b_n^{i,\cdot}(p\cup p_\rev,0)^{1/n},\qq
\beta^{i,\cdot}(p) := \lim_{n\to\oo} b_n^{i\cdot}(p,0)^{1/n}, 
\end{equation*}
exist, and hence the limits exist  in \eqref{eq:324} when $p\notin\Pi$.
Furthermore,
\begin{equation}\label{eq:701}
b_n^{i,\cdot}(p\cup p_\rev,0) \le \beta^{i,\cdot}(p\cup p_\rev)^{n+r},\qq
b_n^{i,\cdot}(p,0) \le \beta^{i,\cdot}(p)^{n+r}.
\end{equation}

\item
If $p\in\Sip^+\setminus\Pi$ and $i\in\{1,2,\dots,M\}$,  
\begin{align*}
\lim_{n\to\oo} b_n^{i,i}(p,0)^{1/n} < \be^+,\qq
\lim_{n\to\oo} b_n^{i,\cdot}(p,0)^{1/n}< \be^+.
\end{align*}
Hence, $\be^+(p)<\be^+(G)$.
\end{letlist}
\end{lemma}

\begin{proof}[Proof of Lemma \ref{lem:10}]
See \cite[Prop.\ 4.2]{GL-loc}. It suffices to consider ($+$)ebridges.
Let 
$B_n=\min\{b_n^{k,\cdot}: k=1,2,\dots,M\}$.
Evidently, $B_{m+n} \ge B_mB_n$, whence the limit $\lambda=\lim_{n\to\oo} B_n^{1/n}$ exists
and, in addition, $B_n \le \lambda^n$.

Let $i\in\{1,2,\dots,M\}$ and pick $j$ such that $b_{n+r}^{j,\cdot} = B_{n+r}$. We extend backwards
a ($+$)ebridge $w$ from $o_i$ with a copy of $\nu(o_j,o_i)$;
more precisely, find $\alpha\in\sH$ such that $\alpha(v_i)=o_i$, and consider the 
ebridge obtained by adjoining $\alpha(\nu(o_j,o_i))$ and $w$. 
On recalling the definition of a \hf, we obtain
\begin{equation}\label{eq:602}
B_n \le b_n^{i,\cdot} 
\le b_{n+l}^{j,\cdot}
\le b_{n+r}^{j,\cdot}
=B_{n+r} \le \lambda^{n+r},
\end{equation}
where $l=\vert\nu(o_j,o_i)\vert$. It follows that, for all $i$,
\begin{equation}\label{eq:601}
(b_n^{i,\cdot})^{1/n} \to \lambda \q\text{as } n\to \oo,
\end{equation}
whence $\beta^+$ exists and satisfies $\beta^+=\lambda$.

For the remaining part of \eqref{eq:603} we argue similarly to \eqref{eq:602}.
By extending an $n$-step $(i,k)$-ebridge 
with a copy of $\nu(o_k,o_j)$, we obtain
\begin{equation*}
b_n^{i,k} \le b_{n+l}^{i,j},
\end{equation*}
where $l=\vert \nu(o_k,o_j)\vert$. On summing over $k$, we find that
\begin{equation*}
b_n^{i,\cdot} \le M \max\{b_m^{i,j}:n \le m\le n+r\}, 
\end{equation*}
and the claim follows.
\end{proof}

\begin{proof}[Proof of Lemma \ref{lem:11}]
(a) This holds as in Lemma \ref{lem:24}, with the same proof 
but subject to the addition of connectors where needed.

(b) Since $p\in\Sip^+$, there exist $a,b\in\{1,2,\dots,M\}$
and an $(a,b)$-ebridge $w$ containing a copy of $p$. By adding suitable connectors at the ends of $w$,
we obtain an $(i,i)$-ebridge containing $p$.

(c) This holds as in the proof of Lemma \ref{lem:10}. 
Follow that proof while restricting to ebridges containing no copy of $p$ or $p_\rev$ (\resp, $p$).
The assumption that $p\notin\Pi$ is used in the step corresponding to \eqref{eq:602}.

(d) Let $p\in\Sip^+\setminus\Pi$ and let 
$q_n^i$ be the number of $n$-step $(i,i)$-ebridges that are $i$-prime. We repeat the argument
leading to  \eqref{eq:110}. First,
\begin{equation*}
b_n^{i,i}=q_n^i+\sum_{m=1}^{n-1} q_m^i b_{n-m}^{i,i}, \qq n\ge 1.
\end{equation*}
Therefore, the generating functions
\begin{equation*}
B^i(x)=\sum_{n=1}^\oo b_n^{i,i} x^{n},\qq Q^i(x)=\sum_{n=1}^\oo q_n^i x^{n},
\end{equation*}
satisfy 
\begin{equation*}
B^i(x) = \frac{Q^i(x)}{1-Q^i(x)}\qq \text{whenever } Q^i(x)<1.
\end{equation*}
By Lemma \ref{lem:10}, $B^i$ has radius of convergence $1/\be^+$, whence
$Q^i(1/\be^+)=1$.

Let $Q^i_p$ be the generating function of $i$-prime $(i,i)$-ebridges that contain no copy
of $p$. 
Since $p$ is assumed ($+$)bridgeable, by part (b),
$Q^i_p (1/\be^{+}) <1$.
Since $Q^i_p$ is continuous inside its ball of convergence, and $M<\oo$, there exists $\eta
=\eta(p)>0$ such that 
\begin{equation}\label{eq:203}
Q^i_p((\be^+)^{-1}+\eta)<1\qq\text{for all } i.
\end{equation}

As following \eqref{eq:103}, we deduce that
\begin{equation*}
\limsup_{n\to\oo} b_n^{i,i}(p,0) ^{1/n} \le \frac1{(\be^+)^{-1}+\eta} < \be^+\qq\text{for all } i.
 \end{equation*}
The limit exists here by subadditivity. This proves the first inequality of part (d).

Moving to the second inequality, 
let $j\in\{1,2,\dots,M\}$. By adding a copy of $\nu(o_j,o_i)$ at the end of
an $(i,j)$-ebridge $w$ of length $m$, we obtain an $(i,i)$-ebridge $w'$ of length $m+l$
where $l=\vert\nu(o_j,o_i)\vert$; moreover, if
$w$ does not contain a copy of $p$, then neither does $w'$
(we use the assumption $p\notin\Pi$ here). Therefore, 
\begin{equation*}
b_n^{i,j}(p,0) \le b_{n+l}^{i,i}(p,0) \le \max\{b_m^{i,i}(p,0): n \le m\le n+r\}. 
\end{equation*}
The second inequality of part (d) follows on summing over $j$ and taking the $n$th root.

With $i$ given by $\id\in\sH o_i$, we have by the above that
\begin{equation*}
\beta^+(p) = \lim_{n\to\oo} b_n^+ (p,0)^{1/n}
=\lim_{n\to\oo} b_n^{i,\cdot}(p,0)^{1/n}
<\beta^+,
\end{equation*}
as claimed.
\end{proof}

\begin{lemma}\label{lem:13}
Let $p\in\Sip^\pm\setminus\Pi$. Then $\mu(G) = \be^+\vee \be^-$, and
\begin{align}\label{eq:new4}
\mu(p\cup p_\rev) &\le \be^+(p\cup p_\rev)\vee\be^-(p\cup p_\rev)\\
 &\le \beta^\pm(p) \vee \be^\mp(p_\rev) < \mu(G).\nonumber
\end{align} 
\end{lemma}

Once this is proved, the proof of the quasi-transitive case of Theorem \ref{thm:pattern1} 
is completed as in the transitive case,
by Lemma \ref{lem:11} followed by an appeal to the proof 
of \cite[Lemma 7.2.5]{ms}. 

\begin{corollary}\label{cor:1}
For $p\in\Sip\setminus\Pi$, the limit 
\begin{equation*}
\mu(p\cup p_\rev) = \lim_{n\to\oo} \sigma_n(p\cup p_\rev, 0)^{1/n}
\end{equation*}
exists and satisfies $\mu(p\cup p_\rev) = \beta^+(p\cup p_\rev)\vee\beta^-(p\cup p_\rev)$.
\end{corollary}

\begin{proof}[Proof of Corollary \ref{cor:1} (using Lemma \ref{lem:13})]
By bounding SAW-counts below by bridge-counts, we have
\begin{equation*}
\sigma_n(p\cup p_\rev) \ge b_n^+(p\cup p_\rev) + b_n^-(p\cup p_\rev),
\end{equation*}
whence 
\begin{equation*}
\liminf_{n\to\oo} \sigma_n(p\cup p_\rev)^{1/n} \ge \beta^+(p\cup p_\rev)\vee\beta^-(p\cup p_\rev).
\end{equation*}
This may be combined with \eqref{eq:new4} to obtain the claim.
\end{proof}

\begin{figure}
\centerline{\includegraphics[width=0.5\textwidth]{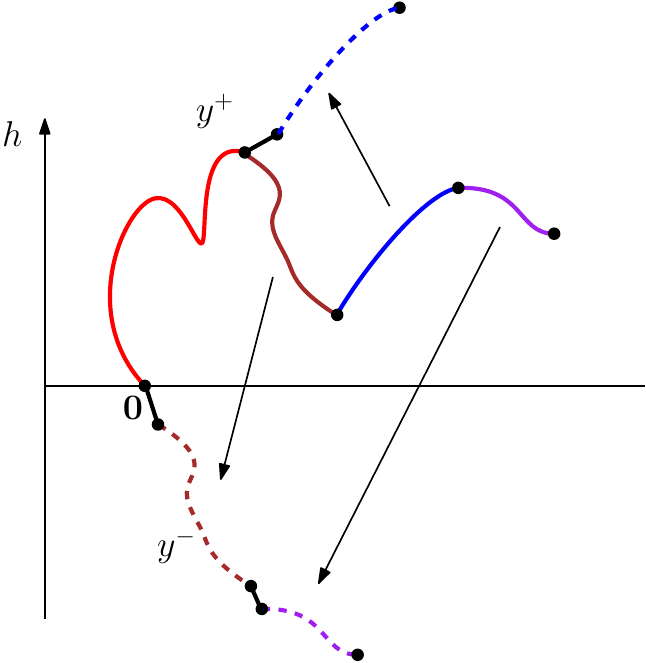}}
\caption{The transformation of Figure \ref{fig:half}
with (black) connectors added where needed.}\label{fig:half2}
\end{figure}

\begin{proof}[Proof of Lemma \ref{lem:13}]
We adapt the arguments of Section \ref{ssec:trans}, noting that, since $p\notin\Pi$, 
no copy of either $p$ or $p_\rev$ appears in any connector. Some further details follow.

In proving \eqref{eq:213}, it was explained how to turn
a half-space walk $w$ (starting at $\id$) into a $(+)$bridge and a $(-)$bridge. This was achieved by 
considering images (under $\sH$) of 
$(+)$sub-bridges and $(-)$sub-bridges of $w$. In the non-transitive case, the endvertices
of the ensuing 
sub-bridges do not generally coincide, and one needs to insert suitable connectors at the joins. This is explained in
\cite{GL-loc, Lin20} to which the reader is referred for further details (see also Figure
\ref{fig:half2}). 
Equation \eqref{eq:212} remains true, but with
\begin{align*}
h_n(a_1,a_2,\dots,a_k) 
&\le (r+1)^{k}
\sum_{m=0}^{n}
\left( \sum_{s=0}^{kr} \sum_{t=0}^\oo  b_{m+s}^+(a_1+a_3+\cdots+t)\right)\\
&\hskip3cm \times\left( \sum_{s=0}^{kr} \sum_{t=0}^\oo b_{n-m+s}^-(a_2+a_4+\cdots+t) \right)\\
&\le (r+1)^{k} \sum_{m=0}^n
\left(\sum_{s=0}^{kr} b^+_{m+s}\right)\left(\sum_{s=0}^{kr} b^-_{n-m+s}\right).
\end{align*}
Therefore, as in \cite[Lemma 3.7]{Lin20} and \eqref{eq:2121}, \eqref{eq:213}, \eqref{eq:bbnd}, 
there exists $B=B(G,h,\sH)<\oo$ such that
\begin{equation*}
h_n \le e^{B\sqrt n} (\be^+\vee\be^-)^n.
\end{equation*}
This bound holds for half-space walks beginning at any given vertex, not only the root.

In a similar manner, \eqref{eq:214} becomes
\begin{align}\label{eq:224}
\si_n &\le \sum_{m=0}^n e^{B\sqrt{n-m}}e^{B\sqrt{m+1}} (\be^+\vee\be^-)^{n+1}\\
&\le e^{C\sqrt n}(\be^+\vee\be^-)^{n} \nonumber
\end{align}
for some $C=C(G,h,\sH)<\oo$.
Therefore, $\mu(G)=\be^+\vee\be^-$.

We may assume without loss of generality that $p\in\Sip^+\setminus\Pi$.
We outline very briefly how to adapt the above sketch in 
order to bound the number of SAWs containing no copy of either
$p$ or $p_\rev$. Corresponding to \eqref{eq:224}, 
we obtain on recalling \eqref{eq:215}--\eqref{eq:new3} and \eqref{eq:701} that
\begin{equation*}
\si_n(p\cup p_\rev,0) \le e^{C\sqrt n}\bigl[\be^+(p\cup p_\rev)\vee\be^-(p\cup p_\rev)\bigr]^n,
 \end{equation*}
whence 
\begin{align*}
\mu(p\cup p_\rev) &\le \be^+(p\cup p_\rev)\vee\be^-(p\cup p_\rev)\\
 &\le \be^+(p)\vee\be^-(p_\rev)< \mu(G),
 \end{align*}
  by 
 \eqref{eq:323} and Lemma \ref{lem:11}(d).
\end{proof}

\subsection{Comments on the proofs of Theorems \ref{thm:t-bridge} and \ref{thm:trail2}}
\label{ssec:comments}

These proofs follow, \emph{mutatis mutandis}, those of Theorems \ref{thm:bridge} and \ref{thm:pattern1},
and they are omitted. 

One difference is noted, namely the following.
Let $v=(v_0,v_1,\dots,v_m)$ and $w=(w_0,w_1,\dots,w_n)$ be SAWs with
$v_m=w_0$, and write $v\circ w$ for the walk obtained by adjoining
$v$ and $w$. In the world of SAWs there are two possibilities, either $v\circ w$ is a SAW or not.
When trails enter into consideration, one must allow the possibility that $v\circ w$ is a trail but not a SAW.
More precisely, copies of a given t-pattern can be created by the act of joining walks together. 
This must be taken into consideration in the surgical aspects of the proof of Theorem \ref{thm:pattern1},
such as in the illustrations of Figures \ref{fig:half}, \ref{fig:full}, and \ref{fig:half2}.

Assume for simplicity that the \hf\ is transitive.
The trail equivalent of Proposition \ref{lem:21}(a) holds in similar fashion. Furthermore, the inequalities
\eqref{eq:2121}--\eqref{eq:215} have equivalents for trails precisely because they are 
\emph{upper bounds} for certain trail-counts.

\section{Theorem \ref{thm:trail}: proof and extension}\label{sec:trailthm}

Let $C$ be a cycle of $G\in\sG$ and $v$ a vertex of $C$. The t-pattern $p(C,v)$ is defined as a trail that
starts and ends at $v$ and traverses $C$ (and no other edge) in between (see Figure \ref{fig:pdelta}).  
(Actually there are two such t-patterns depending on the direction of the traversal, but this
will not be relevant in the following.)
Since $C$ is $2$-edge-connected, $p(C,v)$ is indeed a t-pattern. By \eqref{eq:sym}, 
$\Pi$ comprises SAWs only, so that 
\begin{equation}\label{eq:notinpi}
p(C,v), p(C,v)_\rev \notin \Pi.
\end{equation}

For a graph $G\in\sG$ with a \hf, we present next a necessary and sufficient condition for the inequality  
$\mu<\mut$.
Theorem \ref{thm:trail} is an immediate consequence via the following Proposition \ref{prop:2}. 

Recall the set $\Sitb$ of t-bridgeable t-patterns.

\begin{theorem}\label{prop:3}
Let $d\ge 4$, and let $G\in\sG_d$ possess a \hf\ $(h,\sH)$.
We have that  $\mu<\mut$ if and only if
\begin{equation}\label{eq:criterion}
\text{there exist a cycle $C$ and $v\in C$ such that  $p(C,v) \in \Sitb$.}
\end{equation}

\end{theorem}

\begin{proposition}\label{prop:2}
Let $d\ge 4$, and let $G\in\sG_d$ possess a \emph{transitive} \hf\ $(h,\sH)$.
If $G$ is not a tree, there 
exists a cycle $C$ and a vertex $v$ in $C$, such that $p(C,v) \in \Sitb$.
\end{proposition}

The following question is left open.

\begin{question}\label{qn:weaken}
Does the conclusion of Proposition \ref{prop:2} hold if one weakens the assumption 
on $(h,\sH)$ from transitivity
to quasi-transitivity?
That is, does \eqref{eq:criterion} necessarily hold when the \hf\ is quasi-transitive?
\end{question}

An affirmative answer to this question would permit the 
removal of the assumption in Theorem \ref{thm:trail} that
the height function be transitive. In practical terms, when $G\in\sG$ has a \emph{quasi}-transitive
height function, the criterion of Theorem \ref{prop:3}  can often be checked by inspection.

\begin{proof}[Proof of Proposition \ref{prop:2}]
We shall show the existence of a cycle $C$ such that the following holds. There exists a 
doubly-infinite trail
$t=(\dots,t_{-1},t_0,t_1,\dots)$ that includes a circumnavigation of $C$ and satisfies
\begin{equation*}
h(t_n) \to \begin{cases} -\oo &\text{as } n\to-\oo,\\
\oo &\text{as } n\to\oo.
\end{cases} 
\end{equation*}
By truncating $t$ and adding connectors where needed, we may obtain the required t-bridge containing $p(C,v)$.

Let $D$ be a cycle of $G$. Let $v\in D$ be a vertex of $D$ with maximal height, and let $w\in D$
have minimal height. Since $G$ is transitive and $|D| \ge 2$, we may assume without loss of generality that
\begin{equation}
v\ne w, \qq h(w) \le h(v).
\end{equation}

By the definition of \hf, there exists an infinite \emph{increasing} walk $\pi_v=(v=v_0,v_1,\dots)$ 
of $G$; similarly, there exists an infinite \emph{decreasing} walk $\nu_w=(w=w_0,w_1,\dots)$ 
of $G$. We may choose $\pi_v$ and $\nu_w$ such that
\begin{equation}\label{eq:199}
\begin{aligned}
h(v_1)-h(v)&=\max\{h(v')-h(v): v'\sim v\},\\
h(w)-h(w_1)&=\max\{h(w)-h(w'): w'\sim w\}.
 \end{aligned}
\end{equation}

Since $d\ge 4$, there exists an edge $\langle v,z\rangle$ 
such that $z$ lies in neither the path $\pi_z$ nor the cycle $D$. Let $\nu_z
=(z=z_0,z_1,\dots)$ be an infinite
path of vertices with strictly decreasing heights, and let $z_r$ be the earliest vertex 
of $\pi_z$ lying in $D$ (if no such vertex exists, we set $r=\oo$).
By \eqref{eq:199}, no vertex of $\pi_z$ between $z=z_0$ and $z_r$ lies in $\pi_v\setminus\{v\}$.

There are three cases.
\begin{letlist}
\item
Suppose $r=\oo$. The path $\nu_z$, reversed, followed by the edge $\langle z,v\rangle$,
followed by a traversal of $D$ (in either direction),
followed by $\pi_v$, contains a ($+$)t-bridge that traverses every edge of $D$. The proof is complete in this case
(with $C:=D$).
\item
Suppose $r<\oo$ and $z_r\ne v$.  Let $z_s$ be the last vertex of $\nu_z$ lying in $D$.
The cycle $D$ is the union of two edge-disjoint sub-arcs 
$a_1$, $a_2$ from $z_r$
to $v$, and we label these in such a way that $z_s\in a_1$. 
The sequence $v,z,z_1,z_2,\dots,z_r$, followed by the sub-arc $a_2$,
 constitutes a cycle denoted $C$. The path $(z_s, z_{s+1}, \dots)$, reversed, 
followed by the sub-arc of $a_1$ from $z_s$ to $v$,  followed by a traversal (in either direction) of $C$, followed by $\pi_v$,
constitutes a doubly-infinite trail that traverses every edge of $C$. The proof is complete in this case.
\item 
Suppose $r<\oo$ and $z_r=v$. There are two sub-cases. 
\begin{numlist}
\item If $r\ge 2$, we let $C$ be the cycle $v, z, z_1,\dots,z_r=v$.
The path $\nu_w$, reversed,
followed by a sub-arc of $D$ from $w$ to $v$, followed by a traversal 
of $C$, followed by $\pi_v$, constitutes a doubly-infinite trail that
traverses $C$. The proof is complete in this case.
\item
If $r =1$,  we have that $h(v)<h(z)$. The rest of this proof is devoted to this situation.
\end{numlist}
\end{letlist}

By the above, we may henceforth restrict our attention to the case when
\begin{equation}\label{eq:201-}
\text{for all neighbours $z$ of $v$ not lying in $D\cup \pi_v$, we have $h(v)<h(z)$}.
\end{equation}
By \eqref{eq:201-} and the the maximality of $h(v)$ on $D$, 
\begin{equation}\label{eq:201}
\text{$v$ has exactly $d-2$ neighbours $z$ such that $h(v)<h(z)$}.
\end{equation}
By reversing the roles of $v$ and $w$, we may assume that
\begin{equation}\label{eq:202}
\text{$w$ has exactly $d-2$ neighbours $z$ such that $h(w)>h(z)$},
\end{equation}
in addition to \eqref{eq:201}.

\emph{Assume now that $(h,\sH)$ is transitive.}
Then \eqref{eq:201}--\eqref{eq:202} hold for \emph{all} vertices.
This is a contradiction if $d \ge 5$, and the proof is complete in this case. 

Assume $d=4$. We shall use the notation (in particular, $v$, $r$) 
prior to case (a) above. By \eqref{eq:202} and the transitivity of the \hf, there are
two possible choices for the first vertex $z_1$ of the path $\nu_z$. We may therefore choose $\nu_z$
in such a way that $r\ne 1$.  We now follow the analysis of cases (a)--(c) above to deduce the claim
in this final case.
\end{proof}

Only in the final two paragraphs of the last proof did we used the transitivity of the \hf.

\begin{proof}[Proof of Theorem \ref{prop:3}]
Let $G$ be as given, and assume first that \eqref{eq:criterion} holds.
We may assume without loss of generality that 
$p(C,v)\in\Sitb^+$, and we recall \eqref{eq:notinpi}. By Theorem \ref{thm:trail2}(a),
all but an exponentially small proportion of $n$-step trails from $\id$ contain either $p$ or $p_\rev$.
No such trail is a SAW, whence  $\mu < \mut$.

Suppose conversely that $\mu<\mut$. 
By Theorem \ref{thm:t-bridge}, $\mu< \betat^+\vee \betat^-$. 
We may assume without loss of generality that $\betat^+\ge \betat^-$, whence $\mu < \betat^+$.
Therefore, for large $n$, there exists an $n$-step t-bridge, denoted $w$,  that is not a SAW. 
Since $w$ is not a SAW, it contains some closed sub-walk. It follows that
there exists an appropriate pair $C$, $v$
such that $w$ contains the t-pattern $p(C,v)$; we choose such $C$, $v$ accordingly.
We have that $p(C,v)\in \Sitb$ as required.
\end{proof}

\begin{proof}[Proof of Theorem \ref{thm:trail}]
By Proposition \ref{prop:2}, there exist $C$, $v$ such that $p(C,v)\in\Sitb$. 
The required conclusion follows by Theorem \ref{prop:3}. 
\end{proof}

\section{SAWs on line graphs}\label{sec:line}

\begin{figure}
\centerline{\includegraphics[width=0.5\textwidth]{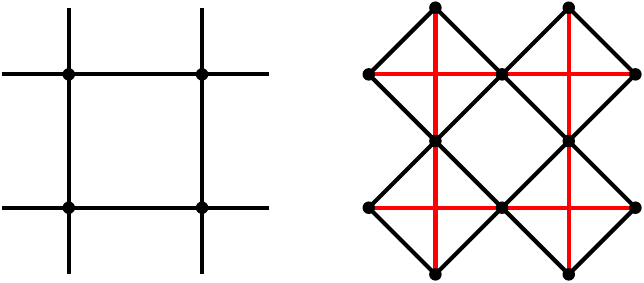}}
\caption{(L) Part of the square lattice $\ZZ^2$, and (R) its line graph $\ZZ^2_K$.
Removal of the red edges inside the squares of $\ZZ^2_K$ results in a copy of $\ZZ^2$.}\label{fig:lineg}
\end{figure}

The \emph{line graph} (or \emph{covering graph}) of a graph $G=(V,E)$ is the graph $G_K=(E,F)$
where $F$ is the set of pairs of elements of $E$ that share a vertex when viewed as edges of $G$. 
We picture the pair $G$, $G_K$ as follows. The vertices of $G_K$ are the midpoints of the edges of $G$, 
and two such vertices are connected in $G_K$ by an edge whenever the corresponding edges of $G$
share a vertex. Thus a vertex $v$ of $G$ with degree $d$ is
replaced by a copy, denoted $K(v)$, of the complete graph $K_d$ on $d$ vertices, having its vertices at the
midpoints of edges of $G$ incident with $v$.   We shall refer to such copies of 
complete graphs as \emph{kites} of $G_K$, and we write $\sK$ for the set of such kites. 
The construction is illustrated in Figure \ref{fig:lineg} in the special case $G=\ZZ^2$.

Let $G\in\sG$. Since $G$ is transitive,  
$G_K$ is quasi-transitive (as in the discussion following \cite[Prop.\ 7.3]{GL20}).
Fix a root $\id$ for $G$, and a root $\id_K \in K(\id)$ for $G_K$.

We introduce three sets of SAWs on $G_K$ starting at $\id_K$. These sets differ in the manner in which they may behave on kites.
\begin{romlist}
\item 
$\Si$ is the set of all SAWs from $\id_K$ on $G_K$.

\item $\Si^1$ is the set of SAWs from $\id_K$ on $G_K$ with the property that, each time they arrive in any given kite $K$,
they must traverse exactly one edge of $K$ before exiting $K$. They may return to $K$ later.

\item
$\Si^-$ is the set of SAWs from $\id_K$ on $G_K$ that may arrive (or begin) 
in any given kite at most once.

\end{romlist}

Write $\mu(G_K)$, $\mu^1(G_K)$, $\mu^-(G_K)$ for the respective exponential growth rates of $n$-step
SAWs in these sets. These constants exist by arguments using subadditivity; see also \cite{jmhII}.

A trail on $G$ gives rise to a SAW on $G_K$ (see \cite[Sect.\ 7.4.2]{Hug}). 
Moreover, the set of such trails generates the set $\Si^1$, whence
\begin{equation}
\mut(G) = \mu^1(G_K).
\end{equation}

We note a simple example of the use of a line graph.

\begin{example}\label{ex:z2k}
Let $G=\ZZ^2$, the square lattice, with line graph $G_K$ as illustrated in Figure \ref{ex:z2k}.
The graph $G_K$ contains a copy $G'$ of $G$.
By \cite[Thm 3.2]{GL14} with $\Ga$ the set of
translations of $G'$, and $\sA$ the subset of elements of $\Ga$ that
preserve $G_K$, we have
that $\mu(G') < \mu(G_K)$.
\end{example}

When restricted to cubic graphs (that is, every vertex-degree is $3$) the above transformation
is the so-called Fisher construction. 
The last was introduced by Michael Fisher \cite{fisher}
in his work on the Ising model, and it has been exploited in \cite{GL13} in the context of
SAWs (see also \cite{GL27,MW}). It was shown in \cite[Thm 1]{GL13} 
that, when $G$ is cubic and quasi-transitive, 
the connective constants of $G$ and $G_K$ are related by a simple formula.
The Fisher construction was adapted in \cite[Sect.\ 7]{GL20}
to transitive graphs with degree $4$, and the ensuing observations were used in \cite[Thm 3.1]{GG26b}
to find upper bounds for connective constants of Grigorchuk graphs. 

Note that
\begin{equation}\label{eq:three}
\mu^-(G_K) =\mu(G_K) \qq\text{when $G$ is cubic}.
\end{equation} 

\begin{theorem}\label{thm:GK}
Let $G\in\sG_d$ where $d\ge 3$.
The connective constants $\mu(G)$ and $\mu^-(G_K)$ satisfy
\begin{equation}\label{eq:triv}
\mu(G_K) \ge \mu^-(G_K), \q\text{and}\q  g_d\left(\frac1{\mu^-(G_K)}\right)= \frac1{\mu(G)},
\end{equation}
where
\begin{equation}\label{eq:triv2}
g_d(\xi) = \sum_{l=1}^{d-1} \xi^l \frac{(d-2)!}{(d-l-1)!},
\end{equation}
with the convention that $0!=1$.
\end{theorem}

In the case $d=3$, equation \eqref{eq:triv2} differs from 
that of \cite[Thm 1(a)]{GL13} owing to the convention 
adopted here that the vertices
of $G_K$ lie at the midpoints of the corresponding edges of $G$.

\begin{proof}
The first inequality of \eqref{eq:triv} is an immediate consequence of  the definitions. 
For the (second) equality, we follow the argument of \cite{GL13}
(see also the proof of \cite[Prop.\ 7.4]{GL20}), 
noting that there are $(d-2)(d-3)\cdots(d-l)$ paths of $K_d$ of length $l$ ($\le d-1$)
joining any two given distinct vertices of $K_d$.
\end{proof}

The relevance of Theorem \ref{thm:GK} to the analysis \cite{GG26b} of SAWs on Grigorchuk graphs is as follows.
Let $G_\om$ be the Cayley graph of the general Grigorchuk group $\Ga_\om$
corresponding to a given sequence $\om\in\{0,1,2\}^\oo$ (see
\cite{GG26b} for definitions). Let $H$ denote the
graph obtained from $G_\om$ by contracting each of its kites (copies of $K_4$ in this case) 
to a single vertex; thus, $G_\om=H_K$ with $d=4$.
In the proof of the main Theorem 3.1 of \cite{GG26b}, a set $W_n$ of 
$n$-step SAWs on $G_\om$ is constructed that is based on its orbital Schreier graph $\sS_\om$.
As noted in \cite[Lemma 5.1]{GG26b}, $W_n\subseteq \Si^-(G_\om)$, whence its growth rate is bounded above
by $\mu^-(G_\om)$ which, in turn, is given by Theorem \ref{thm:GK} in terms of $\mu(H)$.

\section*{Acknowledgement}
The author thanks Gordon Slade for his suggestions, and Stuart Whittington for
sharing his reminiscences of lunchtime conversations with John Hammersley in 1984/5. 

\section*{Declaration}
The author has no conflicts of interest, no data arose in this work, and no use was made of AI.

\providecommand{\bysame}{\leavevmode\hbox to3em{\hrulefill}\thinspace}
\providecommand{\MR}{\relax\ifhmode\unskip\space\fi MR }
\providecommand{\MRhref}[2]{%
  \href{http://www.ams.org/mathscinet-getitem?mr=#1}{#2}
}
\providecommand{\href}[2]{#2}

\end{document}